\documentstyle[11pt]{article}

\font\tinybbfont=msbm6
\font\scriptsizebbfont=msbm6 scaled \magstep 1
 1
\font\footnotesizebbfont=msbm9 scaled \magstep 0
\font\smallbbfont=msbm7 scaled \magstep 2
\font\bbfont=msbm9 scaled \magstep1  
 1
 2
 3
 4

\def\tinyBbb#1{\hbox{\tinybbfont #1}}
\def\scriptsizeBbb#1{\hbox{\scriptsizebbfont #1}}

\def\footnotesizeBbb#1{\hbox{\footnotesizebbfont #1}}
\def\smallBbb#1{\hbox{\smallbbfont #1}}
\def\Bbb#1{\hbox{\bbfont #1}}

 1
 1
 0
 2
 1
 2
 3
 4

\newcommand{\AM}{\mbox{${\cal A}$${\cal M}$}}
\newcommand{\Bun}{\mbox{\it Bun}\,}
\newcommand{\CP}{\mbox{{\Bbb C}{\rm P}}}

\newcommand{\Div}{\mbox{\it Div}}
\newcommand{\Diag}{\mbox{\it Diag}\,}
\newcommand{\Ext}{\mbox{\it Ext}\,}
\newcommand{\Gr}{\mbox{\it Gr}\,}
\newcommand{\Hom}{\mbox{\it Hom}\,}
\newcommand{\Hombf}{\mbox{\bf Hom}\,}
\newcommand{\Id}{\mbox{\it Id}}
\newcommand{\Isom}{\mbox{\it Isom}\,}
\newcommand{\Isombf}{\mbox{\bf Isom}\,}
\newcommand{\Ker}{\mbox{\it Ker}\,}
\newcommand{\Pic}{\mbox{\rm Pic}\,}

\newcommand{\Proj}{\mbox{\rm Proj}\,}
\newcommand{\QM}{\mbox{$\cal Q$${\cal M}$}}
\newcommand{\Sch}{\mbox{\it Sch}\,}
\newcommand{\Sets}{\,\mbox{\it Sets}\,}
\newcommand{\SL}{\mbox{\it SL}}
\newcommand{\SO}{\mbox{\it SO}\,}
\newcommand{\Spec}{\mbox{\it Spec}\,}
\newcommand{\Spin}{\mbox{\it Spin}\,}
\newcommand{\Sym}{\mbox{\it Sym}}
\newcommand{\Tor}{\mbox{\it Tor}\,}
\newcommand{\Tot}{\mbox{\it Tot}\,}
\newcommand{\WD}{\mbox{${\cal W}$${\cal D}$}}

\newcommand{\determinant}{\mbox{\it det}\,}
\newcommand{\dimm}{\mbox{\it dim}\,}
\newcommand{\pre}{\mbox{\it pre}}
\newcommand{\pt}{\mbox{\it pt}}
\newcommand{\rank}{\mbox{\it rank}\,}
\newcommand{\reldeg}{\mbox{\it rel-deg}\,}

\begin{document}

\enlargethispage{23cm}

\begin{titlepage}

$ $

\vspace{-2cm} 

\noindent\hspace{-1cm}
\parbox{6cm}{\small November 2002}\
   \hspace{6.5cm}\
   \parbox{5cm}{math.AG/0212316}

\vspace{3em}

\centerline{\large\bf
 On A-twisted moduli stack for curves from
 }
 \vspace{1ex}
 \centerline{\large\bf
  Witten's gauged linear sigma models
 } 

\vspace{3em}

\hspace{-6em}
\begin{minipage}{18cm}
 \begin{center}
  \parbox[t]{5cm}{
   \centerline{\large Chien-Hao Liu$^1$}
   \vspace{1.1em}
   \centerline{\it Department of Mathematics}
   \centerline{\it Harvard University}
   \centerline{\it Cambridge, MA 02138}  } \
   \hspace{1em}
  \parbox[t]{5cm}{
   \centerline{\large Kefeng Liu$^2$}
   \vspace{1.1em}
   \centerline{\it Department of Mathematics}
   \centerline{\it University of California at Los Angelas}
   \centerline{\it Los Angelas, CA 90095}  } \
   \hspace{1em}
  \parbox[t]{5cm}{
   \centerline{\large Shing-Tung Yau$^3$}
   \vspace{1.1em}
   \centerline{\it Department of Mathematics}
   \centerline{\it Harvard University}
   \centerline{\it Cambridge, MA 02138}  }
 \end{center}
\end{minipage}

\vspace{2em}

\begin{quotation}
\centerline{\bf Abstract}
\vspace{0.3cm}
\baselineskip 12pt  
{\small
 Witten's gauged linear sigma model [Wi1] is one of the universal
  frameworks or structures that lie behind stringy dualities.
 Its A-twisted moduli space at genus $0$ case has been used
  in the Mirror Principle [L-L-Y] that relates Gromov-Witten
  invariants and mirror symmetry computations.
 In this paper the A-twisted moduli stack for higher genus curves
  is defined and systemically studied. It is proved that such
  a moduli stack is an Artin stack.
 For genus $0$, it has the A-twisted moduli space of [M-P] as
  the coarse moduli space.
 The detailed proof of the regularity of the collapsing morphism
  by Jun Li in [L-L-Y$\,$: I and II] can be viewed as a natural morphism
  from the moduli stack of genus $0$ stable maps to
  the A-twisted moduli stack at genus $0$.

 Due to the technical demand of stacks to physicists and the
  conceptual demand of supersymmetry to mathematicians,
  a brief introduction of each topic that is most relevant to the main
  contents of this paper is given in the beginning and the appendix
  respectively.
 Themes for further study are listed in the end.
} 
\end{quotation}

\vspace{1em}

\baselineskip 12pt
{\footnotesize
 } 

\medskip

\noindent {\small
MSC number 2000$\,$: 14D20, 81T30; 14H60, 14M25, 14J32.
} 

\medskip

\baselineskip 11pt
{\footnotesize
\noindent{\bf Acknowledgements.}
 We would like to thank
    David Cox
  for explanations of his work;
    Jun Li
  for discussions of his proof of regularity of the collapsing morphism;
    Bong H.\ Lian
  for the strong influence on our understanding of mirror principle;
    Dan Edidin, Shinobu Hosono, Yi Hu, Daniel Huybrechts, and Chiu-Chu Liu
  for communications on their works;
    Sarah Dean, Kentaro Hori, Shiraz Minwalla, Andrew Strominger,
    and Cumrun Vafa
  for the courses/discussions on the related string theory, 
 Jiun-Cheng Chen, Daniel Freed, J.L., Mihnea Popa
  for comments that lead to corrections and improvements of the draft.
 C.-H.L.\ would like to thank in addition
    Izzet Coskun, Joe Harris, Mircea Mustata, and M.P.\
  for courses/discussions on topics of algebraic geometry;
    Alexander Braverman, J.-C.C., Ian Morrison, Jason Starr,
    Ravi Vakil, and participants of the Workshop on Stacks at MSRI
  for courses/discussions on stacks,
    Orlando Alvarez, Jacques Distler, D.F., and Rafael Nepomechie
  for courses/educations on string/supersymmetry,
  Department of Mathematics at UCLA for hospitalities, and
  Ling-Miao Chou for tremendous moral support.
 The work is supported by NSF grants.
} 

\noindent
\underline{\hspace{20em}}

\noindent
{\footnotesize E-mail:
$^1$chienliu@math.harvard.edu$\,$,\hspace{1ex}
$^2$liu@math.ucla.edu$\,$,\hspace{1ex}
$^3$yau@math.harvard.edu}

\end{titlepage}

%
%

\newpage
$ $

\vspace{-4em}  

\centerline{\sc A-Twisted Moduli Stack from Witten's GLSM}

\vspace{1.5em}

\baselineskip 14pt  

\begin{flushleft}
{\Large\bf 0. Introduction and outline.}
\end{flushleft}

\begin{flushleft}
{\bf Introduction.}
\end{flushleft}
 Witten's gauged linear sigma model (GLSM) [Wi1] is one of
  the universal frameworks or structures that lie behind stringy
  dualities.
 There are many geometric data that are encoded in a GLSM,  
  in particular a toric variety $X$.
 From a gauged linear sigma model, one can obtain two different 
  field theories that have the same local but different global 
  field-theoretic contents as the original theory.
 These descendant theories are called the A-twist and the B-twist of 
  the original theory. 
 The moduli space of vacua of the A-twisted theory is given by 
  solutions to a system of vortex-type equations. 
 Geometrically each solution corresponds to a system of 
  line bundles-with-a-section on curves with these sections
  satisfying some nondegeneracy or nonvanishing conditions.
 (The background from field and string theory are summarized
  in Appendix for mathematicians.)    
  
 In general a system of line bundles-with-a-section has nontrivial
  automorphisms.
 Thus the correct language to study the related moduli problem 
  in the algebro-geometric setting is stack. 
 Since the stack $\AM_g(X)$ in this moduli problem on one hand is 
  related to curves and on the other hand arises from the A-twisted
  theory of a gauged linear sigma model, we will call it the 
  {\it A-twisted moduli stack for curves} (from Witten's GLSM).  

 To really do geometry on a stack, usually one requires it be 
  an algebraic (Artin or Deligne-Mumford) stack. 
 For such stacks, by passing to a covering system of atlases,
  many (down-to-earth) concepts in algebraic geometry for schemes, 
  notably cycles, intersection of cycles, coherent sheaves, 
  and derived categories, can also be defined and studied.   
 (A literature guide on stacks is given in Sec.\ 1 for physicists.)
 In Sec.\ 2, we spell out the definition of $\AM_g(X)$ and prove that
  it is indeed an Artin stack. 
 Hence $\AM_g(X)$ is an object that one may hope to do geometries
  related to curves.  
 
 From the other end, recall the moduli stack 
  $\overline{\cal M}_{g,n}(X)$ of stable maps studied in Behrend
  [Be1, Be2], Behrend-Manin [B-M], Fulton-Pandharipande [F-P], 
  Li-Tian [L-T1, L-T2] $\ldots$ and many others since Kontsevich
  that are related to Gromov-Witten invariants for algebraic varieties.

 It is Witten's insight [Wi1] and Morrison-Plesser's later further
  push [M-P] with some foundation laid down by Cox [Cox2] that the 
  two moduli stacks $\AM_g(X)$ and $\overline{\cal M}_{g,0}(X)$ 
  should be closely related. 
 In particular $\AM_g(X)$ could be as useful in the computation of 
  Gromov-Witten invariants as $\overline{\cal M}_{g,0}(X)$ itself. 
 
 At the moment the detail of relations between these two stacks has
  been carried out for the genus $0$ case. 
  Indeed $\AM_0(X)$, or more precisely its coarse moduli space,
  has been used in the Mirror Principle [L-L-Y$\,$: I and II] that
  relates Gromov-Witten invariants and mirror symmetry computations.
 The two moduli stacks $\AM_0(X)$ and $\overline{\cal M}_{0,0}(X)$
  are related by the natural morphisms
  $$
   \begin{array}{cl}
            \coprod_d M_d(X)   & \\
      \swarrow\hspace{6ex}\searrow  & \\
     \overline{\cal M}_{0,0}(X)\hspace{2em}\AM_0(X) &,\\  
   \end{array}
  $$
  where
   $\coprod_d M_d(X)$ is the moduli stack of genus $0$ stable maps
    into ${\Bbb P}^1\times X$ with the degree on the
    ${\Bbb P}^1$-component being equal to $1$,
   $\coprod_d M_d(X)\rightarrow \overline{\cal M}_{0,0}(X)$
    is the contracting morphism, and 
   $\coprod_d M_d(X)\rightarrow \AM_0(X)$ is the collapsing morphism, 
    whose regularity was proved by Jun Li.
 This is explained more carefully in Sec.\ 3 and Sec.\ 4.
 
 These notes lay down some foundations for several themes to be 
  reported in the future.

\bigskip

\begin{flushleft}
{\bf Outline.}
\end{flushleft}
{\small
\baselineskip 11pt  
\begin{itemize}
 \item [1.]
   A brief tour and literature-guide on stacks for physicists.

 \vspace{-.6ex}
 \item [2.]
  From Cox functor to Witten's A-twisted moduli stack.
  \vspace{-.6ex}
  \begin{itemize}
   \item [2.1]
    The A-twisted moduli stack $\AM_g(X)$.

   \vspace{-.6ex}
   \item [2.2]
    $\AM_g(X)$ is an Artin stack.
  \end{itemize}

 \vspace{-1ex}
 \item [3.]
  The $g=0$ case.

 \vspace{-.6ex}
 \item [4.]
  The collapsing morphism.

 \vspace{-.6ex}
 \item []\hspace{-3ex}Appendix.
  Witten's gauged linear sigma models for mathematicians.   
\end{itemize}
} 

\bigskip

\bigskip

\section{A brief tour and literature-guide on stacks for physicists.}

Basic definitions on stacks needed for the discussions are collected
in this section for introduction of notations and physicists'
convenience.

\bigskip 

\noindent $\bullet$
{\bf Grothendieck topology and site.}
 Let $(\Sch/S_0)$ be the category of schemes over a base scheme $S_0$
  ([Ha1]).
 See [G-M$\,$: Sec.\ II.4], [Kr], and [L-MB$\,$: Chapter 9]
  for the definition - and [Mu1] for why they are needed -
  of the following$\,$:

  {\small
  \begin{itemize}
   \item [$\circ$]
    {\it Topology} = covering system.
    \begin{itemize}
     \item [--]
      {\it \'{e}tale} topology.

     \item [--]
      {\it fppf} topology;
       fppf = faithfully flat (=flat+surjective),
                            + locally of finite presentation.
     \item [--]
     fpqc topology; fpqc = faithfully flat + quasi-compact
    \end{itemize}

   \item [$\circ$]
    {\it Site} on $(Sch/S_0)$
      = usual $(\Sch/S_0)$ + covering systems.
 \end{itemize}
 } 

\bigskip

\noindent
{\it Notation.}
 Let $f:U^{\prime}\rightarrow U$ be a covering of $U$ in the site
  $(\Sch/S_0)$. We shall adopt the following notations for
  the projection maps from fibered products$\,$:
  $$
   U^{\prime\prime\prime}
    :=U^{\prime}\times_U U^{\prime}\times_U U^{\prime}\;
    \stackrel{\pi_{12},\; \pi_{13},\;
     \pi_{23}}{\mbox{---------------}\!\!\longrightarrow}\;
   U^{\prime\prime}:=U^{\prime}\times_U U^{\prime}\;
    \stackrel{p_1,\;
     p_2}{\mbox{---------}\!\!\longrightarrow}\; U^{\prime}\,.
  $$
Such compact notations are particularly useful for diagram-chasings.

\bigskip

\noindent $\bullet$
{\bf Grothendieck's theory of descent.}
Given a covering morphism $f:U^{\prime}\rightarrow U$ in the site
 $(Sch/S_0)$, {\it Grothendieck's theory of descent} studies
 (1) {\it when and how a geometric object} (e.g.\ a coherent sheaf)
     {\it on $U^{\prime}$ can be descended to a geometric object of
     the same kind on $U$} and
 (2) {\it when and how a morphism between descendable geometric
     objects on $U^{\prime}$ descends to a morphism between the
     descent geometric objects on $U$}.
This is a big generalization of the
 {\it local-to-global constructions} in geometry.
See [Kr$\,$: Lecture 4 and Lecture 5] for an introduction and
 references.

\bigskip

\noindent $\bullet$
{\bf Stacks.}
While
 {\it varieties} contain only closed points (= the usual geometric
  points when the ground field $k$ is ${\Bbb C}$),
 {\it schemes} (e.g.\ [E-H] and [Mu4]) contain also nonclosed points
  to make doing geometry more natural.
{\it Stacks} go one step further to contain ``points" with nontrivial
 automorphisms. A ``space" with such a feature is needed
 to parameterize geometric objects that can have nontrivial
 automorphisms, e.g.\ curves and coherent sheaves.
Assuming the background on algebraic geometry in [Ha1], then
 [Mu1], [D-M], [G\'{o}], and [Ed] (in suggested reading order)
 together give a concrete and solid introduction of algebraic stacks
 and their natural appearance in moduli problems in algebraic geometry;
 [L-MB] gives the final up-to-date polishment.
 See also [Art1], [Art2], [Be2], [Bry], [Gil], [H-M], [Mu3], and [Vi]
 for more details.
Recall that a {\it groupoid} is a category in which all the morphisms
 are isomorphisms.

\bigskip

\noindent
{\bf Definition 1.1 [(pre-)stack].} (Cf.\ [D-M], [G\'{o}], and [Kr].)
{\rm
 A {\it stack ${\cal F}$ over $(\Sch/S_0)$} is a
  {\it category fibered in groupoids}
  $p_{\cal F}:{\cal F}\rightarrow(\Sch/S_0)$
 such that the assignment of the fiber
 ${\cal F}(U):=p_{\cal F}^{-1}(U)$ to each $U\in (\Sch/S_0)$
 is a sheaf of groupoids. I.e.\ it is an assignment of groupoids
 $$
  U\in (\Sch/S_0)\; \longrightarrow\; {\cal F}(U)
 $$
 that satisfies the following sheaf axioms$\,$:
 Let $f:U^{\prime}\rightarrow U$ be a covering of $U$ in the site
 $(\Sch/S_0)$.
 \begin{itemize}
  \item [(1)] ({\it Gluing of morphisms}$\,$)\hspace{1em}
   Let ${\cal E}_1$, ${\cal E}_2$ be objects in ${\cal F}(U)$ and
    $\varphi^{\,\prime}: f^{\ast}{\cal E}_1
                          \rightarrow f^{\ast}{\cal E}_2$
    be a morphism in ${\cal F}(U^{\prime})$ such that there exists
    an isomorphism
    $\tau:p_1^{\ast}f^{\ast}{\cal E}_1
      \rightarrow p_2^{\ast}f^{\ast}{\cal E}_2$
    in ${\cal F}(U^{\prime\prime})$.
   Then there exists a morphism
    $\varphi:{\cal E}_1\rightarrow{\cal E}_2$ in ${\cal F}(U)$
    such that $f^{\ast}{\varphi}={\varphi}^{\,\prime}$.

  \item [(2)] ({\it Monopresheaf}$\,$)\hspace{1em}
   Let ${\cal E}_1$, ${\cal E}_2$ be objects in ${\cal F}(U)$ and
    $\varphi$, $\psi: {\cal E}_1\rightarrow{\cal E}_2$ be morphisms
    in ${\cal F}(U)$ such that $f^{\ast}\varphi=f^{\ast}\psi$.
   Then $\varphi=\psi$.

  \item [(3)] ({\it Gluing of objects})\hspace{1em}
   Let ${\cal E}^{\prime}$ be an object in ${\cal F}(U^{\prime})$ and
    $\tau:p_1^{\ast}{\cal E}^{\prime}
                    \rightarrow p_2^{\ast}{\cal E}^{\prime}$
    be an isomorphism in ${\cal F}(U^{\prime\prime})$ such that
    $\pi_{23}^{\ast}\tau\circ\pi_{12}^{\ast}\tau=\pi_{13}^{\ast}\tau$
    in ${\cal F}(U^{\prime\prime\prime})$.
  Then there exists an object ${\cal E}$ in ${\cal F}(U)$ and an
   isomorphism $\sigma:f^{\ast}{\cal E}\rightarrow {\cal E}^{\prime}$
   in ${\cal F}(U^{\prime})$ such that
   $p_2^{\ast}\sigma=\tau\circ p_1^{\ast}\sigma$.
 \end{itemize}
 If ${\cal F}$ satisfies only (1) and (2), then it is called
  a {\it prestack} over $(\Sch/S_0)$.
 In this case, {\it sheafification} (official term$\,$:
  {\it stackification}) of ${\cal F}$ gives a stack canonically
  associated to a prestack
  (cf.\ [Ha1$\,$: II.\ Proposition-Definition 1.2], [Kr$\,$: Lecture 7],
        and [L-MB$\,$: Lemma 3.2].)
} 

\bigskip

\noindent
{\it Remark 1.1.1 {\rm [}category fibered in groupoids{\rm ]}.}
 ([D-M$\,$: Sec.\ 4] and [Kr$\,$: Lecture 6].)
 Given a category ${\cal F}$ over $(\Sch/S_0)$,
  $p_{\cal F}:{\cal F}\rightarrow(\Sch/S_0)$, it is
  {\it fibered in groupoids} over $(\Sch/S_0)$
  if it satisfies the following morphism-lifting properties$\,$:
  \begin{itemize}
   \item [(1)]
    For any $\varphi:U\rightarrow V$ in $(\Sch/S_0)$ and
    $y\in{\cal F}(V)$ there is a map $f:x\rightarrow y$ in ${\cal F}$
    with $p_{\cal F}(f)=\varphi$.

   \item [(2)]
    Given a diagram
     $\begin{array}{rcl}
        x & & \\[-1ex]
          & \hspace{-.8em}\searrow\!\!\!\mbox{\scriptsize
                                      \raisebox{1.6ex}{$f$}} &   \\[-1ex]
          &             &\hspace{-1em} z \\
          & \hspace{-.8em}\nearrow\!\!\!\mbox{\scriptsize
                                    \raisebox{-.6ex}{$g$}} &   \\[-1.4ex]
        y & &
       \end{array}$
     in ${\cal F}$ with its image
     $\begin{array}{rcl}
        U & & \\[-1ex]
          & \hspace{-.8em}\searrow\!\!\!\mbox{\scriptsize
                              \raisebox{1.6ex}{$\varphi$}} &   \\[-1ex]
          &             &\hspace{-1em} W \\
          & \hspace{-.8em}\nearrow\!\!\!\mbox{\scriptsize
                                 \raisebox{-.6ex}{$\psi$}} &   \\[-1.6ex]
        V & &
       \end{array}$
    in $(\Sch/S_0)\,$.
   Then for all $\chi:U\rightarrow V$ such that $\varphi=\psi\circ\chi\,$,
    there is a unique $h:x\rightarrow y$ such that $f=g\circ h$ and
    $p_{\cal F}(h)=\chi$.
  \end{itemize}
 Such lifting properties in many moduli problems, including the one
  studied in the notes, follow automatically by base-change or
  fibered products. Hence we will omit mentioining them.

\bigskip

\noindent
{\it Remark 1.1.2 {\rm [}algebraic $S_0$-spaces{\rm ]}.}
For technical reasons in algebro-geometric study of moduli problems,
 it is natural to introduce the notion of {\it algebraic $S_0$-spaces}
 and use them, instead of schemes, to define atlas and representability
 of morphisms between stacks, cf.\ [L-MB$\,$: chapters 1 and 10].
Such spaces may be thought of as a collection of \'{e}tale (instead of
 Zariski) local charts for a would-be (generally non-existing) scheme.
To keep things down to earth, we do not adapt this convention
 in the notes.

\bigskip

\noindent $\bullet$
{\bf Morphisms between stacks.}
Let $p_{\cal F}:{\cal F}\rightarrow (\Sch/S_0)$ and
    $p_{\cal G}:{\cal G}\rightarrow (\Sch/S_0)$ be stacks over
    $(\Sch/S_0)$.
A morphism from ${\cal F}$ to ${\cal G}$ is a functor
 $F:{\cal F}\rightarrow {\cal G}$ between the two categories
 such that $p_{\cal G}\circ F =p_{\cal F}$.
Explicitly for moduli stacks, this means that $F$ sends a flat family
 of one class of geometric objects to a flat family of another class
 of geometric objects in a way that commutes with base change.
$F$ is {\it representable} if for all $X\in(\Sch/S_0)$ and morphism
 $x:X\rightarrow {\cal G}$, the fibered product (cf.\ next item)
 ${\cal F}\times_{F,{\cal G},\,x}X$ is also in $(\Sch/S_0)$.
Properties of schemes (e.g.\ {\it proper}, {\it separated},
 {\it smooth}, etc.)\ that are stable under base change and of a local
 nature on the target can be defined for representable morphisms of
 stacks via fibered products with schemes$\,$:
 [D-M$\,$: Sec.\ 4], [G\'{o}$\,$: Sec.\ 2], and
 [L-MB$\,$: Definitions (3.9) and (3.10.1)].

\bigskip

\noindent $\bullet$
{\bf Fibered product.}
[Be2$\,$: Lecture 1, Groupoids], [G\'{o}$\,$: Sec.\ 2.2],
 and [L-MB$\,$: Sec.\ (2.2.2)].
Given two morphisms
 $F:{\cal X}\rightarrow{\cal Z}$ and $G:{\cal Y}\rightarrow{\cal Z}$
 of stacks over $(\Sch/S_0)$, their {\it fibered product}
 ${\cal X}\times_{F,{\cal Z},G}{\cal Y}$ (or denoted
 ${\cal X}\times_{\cal Z}{\cal Y}$ when $F$ and $G$ are clear
 from the text) is defined to be the stack over $(\Sch/S_0)$ with
 \begin{itemize}
  \item []
   {\it Objects}$\,$:\hspace{.2ex} \parbox[t]{32em}{
    {\it Triples} $(X, Y, \alpha)$, where $X\in{\cal X}$,
    $Y\in{\cal Y}$, and $\alpha:F(X)\rightarrow G(Y)$
    is an isomorphism in ${\cal Z}$.                     }

  \item []
   \hspace{-4.2ex}
   {\it Morphisms}$\,$:\hspace{.2ex} \parbox[t]{32em}{
    A morphism from $(X_1,Y_1,\alpha_1)$ to $(X_2,Y_2,\alpha_2)$
    is a {\it pair} $(\varphi_{\cal X},\varphi_{\cal Y})$
    of morphisms $\varphi_{\cal X}:X_1\rightarrow X_2$,
                 $\varphi_{\cal Y}:Y_1\rightarrow Y_2$
    over the same morphism $f:U\rightarrow V$ of schemes
    in $(\Sch/S_0)$ such that the following diagram commutes
    $$
     \begin{array}{cccl}
       F(X_1)  & \stackrel{F(\varphi_{\cal X})}{\longrightarrow}
               & F(X_2)                                           & \\
       \mbox{\scriptsize $\alpha_1$}\downarrow\hspace{2ex}
         & & \hspace{2ex}\downarrow\mbox{\scriptsize $\alpha_2$}  & \\
       G(Y_1)  & \stackrel{G(\varphi_{\cal Y})}{\longrightarrow}
               & G(Y_2) & .
     \end{array}
    $$   }
 \end{itemize}


\noindent $\bullet$
$\underline{\Isom}$ and $\Isombf$.
(Cf.\ [D-M$\,$: Definition (I.10)], [Gro$\,$: Sec.\ 4], and
      [Mu1: Sec.\ 3].)
Given a pair of families of geometric objects, e.g.\ stable curves,
 $\pi_i:X_i\rightarrow S_i$, $i=1,\,2$, then each induces a family
 of geometric object, still denoted by $\pi_i$, over
 $S_1\times_{S_0}S_2$ via pullback.
Let $(\Sets)$ be the category of sets.
Then $\underline{\Isom}(\pi_1,\pi_2)$ is the functor
 $$
  \begin{array}{ccccc}
   \underline{\Isom}(\pi_1,\pi_2) & : & (\Sch/S_0)
                         & \longrightarrow & (\Sets) \\[.6ex]
    &  & S & \longmapsto
       & \{\, (\alpha,\beta)\,|\,
         \alpha\in Hom(S, S_1\times_{S_0}S_2)\,,\;
         \beta: \alpha^{\ast}\pi_1\simeq\alpha^{\ast}\pi_2\,\}\,.
  \end{array}
 $$
In case $\underline{\Isom}(\pi_1,\pi_2)$ is a representable functor,
 the scheme that represents $\underline{\Isom}(\pi_1,\pi_2)$
 will be denoted by $\Isombf(\pi_1,\pi_2)$.
Representability of $\Isom$ in many moduli problems boils down to
 the representability of the Hilbert functor or the Quot functor
 ([Gro]).
When the moduli problem is described by a stack ${\cal X}$ over
 $(\Sch/S_0)$, then $\pi_i$ correspond to morphisms
 $F_i:S\rightarrow {\cal X}$ and
 $\underline{\Isom}(\pi_1,\pi_2)= S_1\times_{F_1,{\cal X},F_2}S_2\,$.
Similarly for $\underline{\Hom}$ and $\Hombf$ that replace
 isomorphisms in the definition of $\underline{\Isom}$ and $\Isombf$
 by morphisms.

\bigskip

\noindent
{\bf Definition 1.2 [Artin stack].}
 ([G\'{o}$\,$: Definition 2.22].) {\rm
 An {\it Artin stack} ${\cal F}$ over $(\Sch/S_0)$ is stack ${\cal F}$
  over $(\Sch/S_0)$ that satisfies additional conditions$\,$:
  \begin{itemize}
   \item [(1)]
    The diagonal morphism
     $\Delta_{\cal F}\rightarrow {\cal F}
        \times_{\mbox{(\scriptsize\it Sch}/S_0)}{\cal F}$
     is representable, quasi-compact, and separated.

   \item [(2)]
    There exists a scheme $U$ - called an {\it atlas} - and
    a smooth and surjective morphism $u:U\rightarrow {\cal F}$.
  \end{itemize}
 See [L-MB$\,$: Definition (5.2)] for the definition of the
  {\it set of points} $|{\cal F}\,|$ of a stack ${\cal F}$.
} 

\bigskip

\noindent
{\bf Definiton 1.3 [coarse moduli space].}
 (Cf.\ [G\'{o}$\,$: Definition 2.6] and [Vi$\,$: (2.1) Definition].)
{\rm
 A {\it coarse moduli space for a stack} ${\cal F}$ is a scheme $Z$
  together with a morphism $\phi: {\cal F}\rightarrow Z$ such that
  \begin{itemize}
   \item [$(i)$]
     if $Z^{\prime}$ is another scheme that admits a morphism
     $\phi^{\prime}:{\cal F}\rightarrow Z^{\prime}$ then there is
     a unique morphism of schemes $\eta:Z\rightarrow Z^{\prime}$
     with $\phi^{\prime}=\eta\circ\phi\,$,
     (i.e.\ $Z$ {\it coreprsents} ${\cal F}$).

   \item [$(ii)$]
     for any algebraically closed field $k$, the induecd map on
     $k$-points $|\phi\,|:|F\,|(k)\rightarrow Z(k)$ is bijective.
  \end{itemize}
 (Thus, when exists, $Z$ is unique up to a canonical isomorphism.)
} 

\bigskip

\noindent $\bullet$
{\bf Quotient stack.}
 [Be2$\,$: Lecture 1, Example 18.3 and Example 20.4 ],
  [G\'{o}$\,$: Example 2.14], and [L-MB$\,$: Sec.\ (2.4.2)].
 For $S_0=\Spec k$, where $k$ is a ground field, let $G$ be an
  algebraic group over $k$.
 The {\it quotient stack} of a $G$-action on a $k$-scheme $X$
  is denoted by $[X/G]$. It is the stackification of the prestack
  $\pre[X/G]$ over $(\Sch/k)$. An object of the groupoid
  $\pre[X/G](U)$, $U\in(\Sch/k)$, is a diagram
  $$
   \begin{array}{ccc}
     P & \stackrel{f}{\longrightarrow}  & X \\
     \downarrow   & &  \\
     U            & & \hspace{1em},
   \end{array}
  $$
  where $P$ is a principal $G$-bundle over $U$ and
  $f$ is a $G$-equivariant $k$-morphism.

\bigskip

\bigskip

\section{From Cox functor to Witten's A-twisted moduli stack.}
Notations and terminologies of toric geometry used here follow
 mainly [Fu], see also [Oda].

\bigskip

\subsection{The A-twisted moduli stack $\AM_g(X)$.}

Let
 $N\simeq{\Bbb Z}^n$ be a lattice,
 $M$ be its dual lattice,
 $\Delta$ be a fan in $N_{\scriptsizeBbb R}$,
 $\Delta(1)$ be the $1$-dimensional cones of $\Delta$, and
 $n_{\rho}$ be the generator of $\rho\cap N$ for $\rho\in\Delta(1)$.
Let $X$ be the smooth toric variety associated to $\Delta$ and $Y$
 be a scheme over $S$.
Recall the following definition from [Cox2]$\,$:

\bigskip

\noindent
{\bf Definition 2.1.1 [$\Delta$-collection].} {\rm
 A {\it $\Delta$-collection} $(L_{\rho},u_{\rho},c_m)_{\rho,m}$
  on $Y/S$ consists of line bundles $L_{\rho}$ on $Y$ flat over $S$,
  sections $u_{\rho}\in H^0(Y,L_{\rho})$ indexed by $\rho\in\Delta(1)$,
  and a collection of isomorphisms
  $c_m:\otimes_{\rho}L_{\rho}^{\otimes\langle m,n_{\rho}\rangle}
                    \simeq {\cal O}_Y$ indexed by $m\in M$ such that
 \begin{itemize}
  \item [$(i)$] \parbox{17ex}{\it Compatibility$\,:$}
   \parbox[t]{10cm}{$c_m\otimes c_{m^{\prime}}=c_{m+m^{\prime}}$
                    for all $m, m^{\prime}\in M$.
             }

  \item [$(ii)$] \parbox{17ex}{\it Nondegeneracy$\,:$}
   \parbox[t]{10cm}{
    The map
     $\;\sum_{\sigma\in\Delta_{max}}\otimes_{\rho\not\subset\sigma}
        u_{\rho}^{\ast}\,:\,
       \oplus_{\sigma\in\Delta_{max}}\otimes_{\rho\not\subset\sigma}
        L_{\rho}^{-1}\, \rightarrow\, {\cal O}_Y\;$
     is surjective, where
    $u_{\rho}^{\ast}:L_{\rho}^{-1}\rightarrow {\cal O}_Y$ is the dual
    morphism of $u_{\rho}:{\cal O}_Y\rightarrow L_{\rho}$.
    }
 \end{itemize}
 An {\it isomorphism}
 $(L_{\rho},u_{\rho},c_m)_{\rho,m}
   \stackrel{\sim}{\rightarrow}
   (L_{\rho}^{\prime}, u_{\rho}^{\prime}, c_m^{\prime})_{\rho,m}$
 consists of isomorphisms
 $\gamma_{\rho}:L_{\rho}\stackrel{\sim}{\rightarrow}L_{\rho}^{\prime}$
 which carry $u_{\rho}$ to $u_{\rho}^{\prime}$ and $c_m$ to $c_m^{\prime}$.
} 

\bigskip

\noindent
{\bf Explanation/Fact 2.1.2 [Cox].} {\rm
 For the application in this article, we will consider only the case
  when the set $\{\,n_{\rho}\,\}_{\rho}$ spans $N_{\scriptsizeBbb R}$.
 In this case
  $X=({\Bbb C}^{\Delta(1)}-V(I))/\mbox{\raisebox{-.4ex}{$G$}}\,$,
  where ${\Bbb C}^{\Delta(1)}=\Spec{\Bbb C}[x_{\rho}:\rho\in\Delta(1)]$,
   $I$ is the ideal generated by
    $\prod_{\rho\not\subset\sigma}\,x_{\rho}$, $\sigma\in\Delta_{max}$,
  and $G=\Hom_{\scriptsizeBbb Z}(\Pic(X),{\Bbb C}^{\times})$ acts on
  ${\Bbb C}^{\Delta((1)}$ via the exact sequence
  $$
   1 \;  \longrightarrow\;  G\;
     \longrightarrow\;
      \Hom_{\scriptsizeBbb Z}({\Bbb Z}^{\Delta(1)},{\Bbb C}^{\times})\;
     \longrightarrow\; T_N \; \longrightarrow\; 1\,.
  $$
 The following statements either are explicitly in or follow
  immediately from {\rm [Cox2]}.
 \begin{itemize}
  \item [(1)]
  {\it The isomorphisms $c_m$}.
    \begin{itemize}
     \item [(1.1)]
      $L_{\rho}$ are unrelated abstract line bundles on $Y/S$.
       To relate the collection of sections $u_{\rho}$ to a map from
       $Y$ to ${\Bbb C}^{\Delta(1)}$, some data is needed that enables
       one to compare sections in different $L_{\rho}$
       - more precisely, the induced sections from $u_{\rho}$
         on isomorphic tensor products of $L_{\rho}$ -.
      The data $c_m$ gives exactly this information up to the
       $G$-action. Condition $(i)$ (Compatibility) is the cocycle
       conditions that make sure this comparison of sections on
       different $L_{\rho}$ is consistent among themselves.

     \medskip
     \item [(1.2)]
      Given $\{L_{\rho}\}_{\rho}$ and two choices $\{c_m\}_m$
       and $\{c_m^{\prime}\}_m$, there exist automorphisms
       $\gamma_{\rho}$ on $L_{\rho}$ that carry $\{c_m\}_m$
       to $\{c_m^{\prime}\}_m$.
      Thus, up to isomorphisms, there is exactly one way to compare
       the line bundles $L_{\rho}$.

     \medskip
     \item []
      {\it Reason.} (Cf.\ [Cox2$\,$: Theorem 1.1, proof].)
       A pair of collections of isomorphisms
        $(\{c_m\}_m, \{c_m^{\prime}\}_m)$ determines a morphism
        $\alpha: M\rightarrow H^0(Y,{\cal O}_Y^{\,\ast})$.
       From the long exact sequence

       \vspace{-1ex}
       \item []
       {\footnotesize
       $$
        \begin{array}{l}
         0\rightarrow \Hom(\Pic X, H^0(Y,{\cal O}_Y^{\,\ast}))
          \rightarrow
           \Hom({\Bbb Z}^{\Delta(1)}, H^0(Y,{\cal O}_Y^{\,\ast}))
          \rightarrow
           \Hom(M, H^0(Y,{\cal O}_Y^{\,\ast})) \\[.6ex]
         \hspace{1em}\rightarrow
           \Ext^1(\Pic X, H^0(Y,{\cal O}_Y^{\,\ast}))\;\; (=0)\;
          \rightarrow \cdots
        \end{array}
       $$
       {\normalsize induced}} 
       from the short exact sequence
        $0\rightarrow M \rightarrow {\Bbb Z}^{\Delta(1)}
          \rightarrow \Pic X \rightarrow 0$,
       one concludes that $\alpha$ can be lifted to a morphism
       $\widetilde{\alpha}: {\Bbb Z}^{\Delta(1)}
         \rightarrow H^0(Y,{\cal O}_Y^{\,\ast})$.
       The morphism $\widetilde{\alpha}$ defines then a collection
       of automorphisms
       $\gamma_{\rho}$ on $L_{\rho}$, $\rho\in\Delta(1)$,
       that carry $\{c_m\}_m$ to $\{c_m^{\prime}\}_m\,$.
    \end{itemize}
 \end{itemize}
 \hspace{14cm}$\Box$

\begin{itemize}
  \item [(2)]
  {\it The nondegeneracy condition}.
   Recall (e.g.\ [Ha1: Appendix A.3]) that, given a section
    $u_{\rho}:{\cal O}_Y\rightarrow L_{\rho}$,
    the zero-sheme of $u_{\rho}$ is defined by the ideal sheaf
    $u_{\rho}^{\ast}(L_{\rho}^{-1})$.
   Thus, Condition $(ii)$ (Nondegeneracy) of a $\Delta$-collection
    means exactly that the image of the map
    $U\rightarrow {\Bbb C}^{\Delta(1)}$ in Item (2) above lies
    completely in ${\Bbb C}^{\Delta(1)}-V(I)$.
 \end{itemize}
} 

\bigskip

Explanation/Fact 2.1.2$\,$: Item (2), [Wi1] and [M-P] together
 lead to the following definitions.

\bigskip

\noindent
{\bf Definition 2.1.3 [weak $\Delta$-collection].} {\rm
 (1)
 A weak $\Delta$-collection on $Y/S$ is a set of data
 $(L_{\rho},u_{\rho},c_m)_{\rho,m}$ as in Definition 2.1.1
 with Condition $(ii)$ (Nondegeneracy) replaced by
 \begin{itemize}
  \item [$(ii^{\,\prime})$] \parbox{17ex}{\it Nonvanishing$\,:$}
   \parbox[t]{10cm}{
    The map
    $\;\sum_{\sigma\in\Delta_{max}}\otimes_{\rho\not\subset\sigma}
       u_{\rho}^{\ast}\,:\,
      \oplus_{\sigma\in\Delta_{max}}\otimes_{\rho\not\subset\sigma}
       L_{\rho}^{-1}\, \rightarrow\, {\cal O}_Y\;$
    is not a zero-morphism when restricted to each irreducible
    component of fibers $Y_s$ over $s\in S$.
    } 
 \end{itemize}
 Isomorphisms of such data are defined the same as in Definition 2.1.1.

 \medskip

 \noindent
 (2) Let $(L_{\rho}, u_{\rho}, c_m)$
      (resp.\ $(L_{\rho}^{\prime},u_{\rho}^{\prime},c_m^{\prime})$)
     be a weak $\Delta$-collection on $Y/S$
      (resp.\ $Y^{\prime}/S^{\prime}$).
     Then a {\it morphism} from $(Y/S, (L_{\rho}, u_{\rho},c_m))$ to
      $(Y^{\prime}/S^{\prime},
          (L_{\rho}^{\prime},u_{\rho}^{\prime},c_m^{\prime}))$
      is a pair $(f,\gamma)$, where $f:Y\rightarrow Y^{\prime}$
      fits into a commutative diagram
     $$
      \begin{array}{ccc}
       Y  & \stackrel{f}{\longrightarrow}  & Y^{\prime} \\
       \downarrow & & \downarrow      \\
       S  & \stackrel{\overline{f}}{\longrightarrow}  & S^{\prime}
      \end{array}
     $$
     and
     $\gamma:
      (L_{\rho}, u_{\rho}, c_m)\stackrel{\sim}{\rightarrow}
       f^{\ast}(L_{\rho}^{\prime},u_{\rho}^{\prime},c_m^{\prime})$
     on $Y/S$.
} 

\bigskip

\noindent
{\bf Definition 2.1.4 [quasistable curves over $S$].} (Cf.\ [Ca].) {\rm
 A prestable (i.e.\ reduced connected nodal) curve is called
  {\it quasistable} if all its destabilizing chains have length $1$.
 A {\it quasistable curve over $S$} is a flat family
  $\pi:{\cal C}\rightarrow S$ of quasistable curves over $S$.
 Define $\QM_g$ to be the category fibered in groupoids of
  quasistable curves over $(\Sch/S_0)$.
} 

\bigskip

\noindent
{\bf Definition/Lemma 2.1.5 [$\AM_g(X)$ stack].} {\it
 Let $(\Sch/S_0)$ be equipped with the fpqc or the fppf topology.
 Define $\AM_g(X)$ to be the category over $(\Sch/S_0)$
  whose fiber over $U\in(\Sch/S_0)$ is given by the groupoid
  $$
   \AM_g(X)(U)\;=\; \{\,
    \mbox{weak $\Delta$-collections on quasistable curves ${\cal C}$
          over $U$}\,\}\,.
  $$
  Then $\AM_g(X)$ is a stack. We shall call it the
  {\rm A-twisted moduli stack associated to $X$ for genus $g$ curves}.
} 

\bigskip

\noindent
{\it Remark 2.1.6.}
 Compared with [Wi1$\,$: Sec.\ 3.4] and [M-P$\,$: Sec.\ 3.7] summarized
  in Appendix, $\AM_g(X)$ is related to the moduli space of
  the A-twisted gauged linear model in the higher genus case.

\bigskip

Definition/Lemma 2.1.5 follows from the proof of the effectiveness
 of a descent datum in the case of quasi-coherent sheaves on schemes
 in $(\Sch/S_0)$, which we recall from [Kr].
See also [SGA1].

\bigskip

\noindent
{\bf Fact 2.1.7 [descent of quasi-coherent sheaves].} {\it
 Let $f:U^{\prime}\rightarrow U$ be a fpqc or fppf morphism in
  $(\Sch/S_0)$. Recall the projection maps from Sec.\ 1
  $$
   U^{\prime\prime\prime}
    :=U^{\prime}\times_U U^{\prime}\times_U U^{\prime}\;
    \stackrel{\pi_{12},\; \pi_{13},\;
     \pi_{23}}{\mbox{---------------}\!\!\longrightarrow}\;
   U^{\prime\prime}:=U^{\prime}\times_U U^{\prime}\;
    \stackrel{p_1,\;
     p_2}{\mbox{---------}\!\!\longrightarrow}\; U^{\prime}\,.
  $$
 \begin{itemize}
  \item [$(a)$]
   {\rm Descent of quasi-coherent sheaves.}
    Let
     ${\cal E}^{\prime}$ be a quasi-coherent
       ${\cal O}_{U^{\prime}}$-module and
     $\tau:p_1^{\ast}\,{\cal E}^{\prime}
           \rightarrow p_2^{\ast}\,{\cal E}^{\prime}$ be an isomorphism
     that satisfies
     $\pi_{23}^{\ast}\tau\circ\pi_{12}^{\ast}\tau =\pi_{13}^{\ast}\tau$.
    Then there exists a quasi-coherent ${\cal O}_U$-module ${\cal E}$
     on $U$ together with an isomorphism
     $\sigma:f^{\ast}\,{\cal E}\rightarrow{\cal E}^{\prime}$
     such that $p_2^{\ast}\sigma=\tau\circ p_1^{\ast}\sigma$.
    The sheaf ${\cal E}$ is unique up to a canonical isomorphism.

  \item [$(b)$]
   {\rm Descent of morphisms.}
   Let
    $({\cal E}^{\prime},\tau)$ and $({\cal F}^{\prime},\upsilon)$
     be descent data and
    $({\cal E},\sigma)$ and $({\cal F},\rho)$ be their respective
     descent as in Item {\rm ({\it a})}.
   Let
    $h^{\prime}:{\cal E}^{\prime}\rightarrow{\cal F}^{\prime}$
     be a morphism that satisfies
     $p_2^{\ast}h^{\prime}\circ\tau=\upsilon\circ p_1^{\ast}h^{\prime}$.
   Then there exists a unique morphism $h:{\cal E}\rightarrow{\cal F}$
    such that $\rho\circ f^{\ast}h=h^{\prime}\circ\sigma$.
 \end{itemize}
} 

\bigskip

\noindent
{\it Sketch of proof.}
Consider the following stricter version of Statement ({\it b})$\,$:
\begin{itemize}
 \item [$(b^{\ast})$]
  {\it Let
        ${\cal E}$ and ${\cal F}$ be quasi-coherent sheaves on $U$ and
        $h^{\prime}:f^{\ast}\,{\cal E}\rightarrow f^{\ast}\,{\cal F}$
        be a morphism of ${\cal O}_{U^{\prime}}$-modules such that
        $p_1^{\ast}h^{\prime}=p_2^{\ast}h^{\prime}$.
       Then there exists a unique morphism $h:{\cal E}\rightarrow{\cal F}$
        such that $f^{\ast}h=h^{\prime}$.     }
\end{itemize}
Via diagram chasings, Statement $(b)$ follows from
 Statement $(b^{\ast})$ and the existence part of Statement $(a)$.
The proof now consists of three steps.

\bigskip

\noindent
(1) {\it Case}$\,$:
    {\it  $f\,=$ faithfully flat morphism between affine schemes.}

\medskip

\noindent
$(1.a)$ {\it Descent of quasi-coherent sheaves.}

\medskip

\noindent
 Define
  $\overline{f}:= f\circ p_1$ ($=f\circ p_2$) and let
  $p_1^{\,\sharp}$,
  $p_2^{\,\sharp}:{\cal O}_{U^{\prime}}
                 \rightarrow {\cal O}_{U^{\prime\prime}}$
  be the defining ring morphism of structure sheaves associated
  to $p_1$, $p_2$ respectively.
 Then their difference $p_1^{\,\sharp}-p_2^{\,\sharp}$
  defines an ${\cal O}_U$-module morphism
  $q:f_{\ast}{\cal O}_{U^{\prime}}
     \rightarrow \overline{f}_{\ast}{\cal O}_{U^{\prime\prime}}$.
 These morphisms induce natural morphisms
  $\widehat{p_1}:{\cal E}^{\prime}
           \rightarrow p_1^{\ast}\,{\cal E}^{\prime}$
  and
  $\widehat{p_2}:{\cal E}^{\prime}
           \rightarrow p_2^{\ast}\,{\cal E}^{\prime}$
  by pulling back global sections.

 Let
   $\Theta := \tau\circ\widehat{p_1}-\widehat{p_2}:
     {\cal E}^{\prime}\rightarrow p_2^{\ast}\,{\cal E}^{\prime}$.
  Then the descent ${\cal E}$ of the descent datum
  $({\cal E}^{\prime},\tau)$ is given by
  ${\cal E}=f_{\ast}\,\Ker\Theta$.
  By definition ${\cal E}$ fits into the following exact sequence
  of ${\cal O}_U$-modules
  $$
   0\; \longrightarrow\; {\cal E}\;
       \stackrel{\iota}{\longrightarrow}\;f_{\ast}\,{\cal E}^{\prime}\;
       \stackrel{\Theta}{\longrightarrow}\;
         \overline{f}_{\ast}\,(p_2^{\ast}\,{\cal E}^{\prime})\,.
  $$
 $\tau$ determines
  an isomorphism
  $k:f^{\ast}\,\overline{f}_{\ast}\,p_2^{\ast}\,{\cal E}^{\prime}\,
    \rightarrow\,
    {p_2}_{\ast}\,\overline{f}^{\ast}\,f_{\ast}\,{\cal E}^{\prime}$
  and an automorphism $h$ on $f^{\ast}\,f_{\ast}\,{\cal E}^{\prime}$.
 These fit into a commutative square of ${\cal O}_{U^{\prime}}$-modules
  $$
   \begin{array}{cccccccl}
    0 & \longrightarrow  & f^{\ast}\,{\cal E}
      & \stackrel{f^{\ast}\iota}{\longrightarrow}
      & f^{\ast}\,f_{\ast}\,{\cal E}^{\prime}
      & \stackrel{f^{\ast}\Theta}{\longrightarrow}
      & f^{\ast}\,\overline{f}_{\ast}\,p_2^{\ast}\,{\cal E}^{\prime}
      & \\
    & & & & \downarrow\,\mbox{\scriptsize $h$} &
      & \downarrow\,\mbox{\scriptsize $k$}     &  \\[-1ex]
    0 & \longrightarrow  & {\cal E}^{\prime}
      & \stackrel{\widetilde{f}}{\longrightarrow}
      & f^{\ast}\,f_{\ast}\,{\cal E}^{\prime}
      & \stackrel{\widetilde{q}}{\longrightarrow}
      & {p_2}_{\ast}\,\overline{f}^{\ast}\,f_{\ast}\,{\cal E}^{\prime}
      & ,
   \end{array}
  $$
  where
   $\widetilde{f}$ and $\widetilde{q}$ are natural morphisms induced
    by $f$ and $q$ respectively and
   both horizontal complexes are exact.
 This determines an isomorphism
  $\sigma:f^{\ast}\,{\cal E}\rightarrow{\cal E}^{\prime}$
  that satisfies the conditions in Statement (a).
  It has the property that if $s^{\prime}$ is a global section in
  $\Ker\,\Theta$ and $s$ is its corresponding global section in
  ${\cal E}$, then $\sigma(f^{\ast}s)=s^{\prime}$.
 Uniqueness of ${\cal E}$ up to a canonical isomorphism follows
  from Item (1.b$^{\ast}$) below.

\bigskip

\noindent
$(1.b^{\ast})$ {\it Descent of morphisms.}

\medskip

\noindent
 The following natural sequence of ${\cal O}_U$-modules induced
  by $f$ and $q$ is exact
  $$
    0\; \longrightarrow\; {\cal F}\;
        \stackrel{\widehat{f}}\longrightarrow\;
           f_{\ast}\,f^{\ast}\,{\cal F}\;
        \stackrel{\widehat{q}}\longrightarrow\;
         \overline{f}_{\ast}\,\overline{f}^{\ast}\,{\cal F}\,.
  $$
 Hence for any global section $s$ in ${\cal E}$, let $f^{\ast}(s)$
  be the corresponding global section in $f^{\ast}\,{\cal E}$,
  then $h^{\prime}(f^{\ast}(s))=f^{\ast}(t)$ for a unique global section
  in ${\cal F}$.

 One can now define the descent morphism $h:{\cal E}\rightarrow{\cal F}$
  by setting $h(s)=t$. By definition $f^{\ast}h=h^{\prime}$.
 Since $f^{\ast}$ is an exact and faithful functor,
  such $h$ is unique.

\bigskip

\noindent
(2) {\it Statements for reductions.}

\medskip

\noindent
 Let
  $R\stackrel{f_1}{\longrightarrow}S\stackrel{f_2}{\longrightarrow}T$
  be a chain of morphisms of schemes in $(\Sch/S_0)$.
 By chasing the following two diagrams$\,$:
  $$
   \begin{array}{ccccccc}
    R\times_S R\times_S R  & \longrightarrow &  R\times_T R\times_T R
      & \longrightarrow & S\times_T S\times_T S  & & \\
    \downarrow\!\downarrow\!\downarrow &
      & \downarrow\!\downarrow\!\downarrow
      & & \downarrow\!\downarrow\!\downarrow &&          \\
    R\times_S R & \longrightarrow & R\times_T R
      & \longrightarrow & S\times_T S  & &           \\
    & \searrow\hspace{-1.6ex}\searrow  & \downarrow\!\downarrow &
     & \downarrow\!\downarrow & &                      \\
    & & R & \stackrel{f_1}{\longrightarrow} & S
          & \stackrel{f_2}{\longrightarrow} & T
   \end{array}
  $$
 and
  $$
   \begin{array}{ccc}
    (R\times_S R)\times_T(R\times_S R)
      & \simeq & (R\times_T R)\times_{S\times_T S}(R\times_T R) \\
    \downarrow\!\downarrow & & \downarrow\!\downarrow \\
    R\times_S R & \longrightarrow & R\times_T R\,,
   \end{array}
  $$
 %
 %
 %
 where all the morphisms are natural projection maps from
 fibered products, one concludes the following statements
 for reduction$\,$:

\bigskip

\noindent
 {\bf Reduction (2.1).} {\it
 Suppose that $(b^{\ast})$ holds for $f_1$ and that for any
  quasi-coherent sheaves ${\cal A}$, ${\cal B}$ on $S\times_T S$,
  the map
  $\eta^{\ast}:\Hom_{S\times_T S}({\cal A},{\cal B})\rightarrow
   \Hom_{R\times_T R}(\eta^{\ast}{\cal A},\eta^{\ast}{\cal B})$
  is injective, then $(b^{\ast})$ holds for $f_2$
  if and only if $(b^{\ast})$ holds for $f_2\circ f_1$.    }

 \bigskip

 \noindent
 {\bf Reduction (2.2).} {\it
 Suppose that $(b^{\ast})$ hold for $f_1$ as well as any pull-back
  of $f_1$. Suppose also that $(a)$ and $(b)$ hold for $f_1$,
  then $(a)$ and $(b)$ hold for $f_2$ if and only if $(a)$ and $(b)$
  hold for $f_2\circ f_1$.  }

\bigskip

\noindent
(3) {\it General case by reductions.}

\medskip

\noindent
 For a general fpqc or fppf morphism $f:U^{\prime}\rightarrow U$
  between $S_0$-schemes, let $\{\,V_i\,\}_{i\in I}$ be a Zariski affine
  cover of $U$ and for each $i\in I$ let
  $\{\,V_{i,j}^{\prime}\,\}_{j\in J_i}$
  be a Zariski affine cover of $f^{-1}(V_i)$.
 If $f$ is fpqc, then $J_i$ can be assumed to be a finite set.
  Then for each $i$ the map
  $f_i:\coprod_j V_{i,j}^{\prime}\rightarrow V_i$
  is a faithfully flat affine morphism and hence $(a)$ and $(b)$ hold
  for $f_i$. Applying Reduction (2.1) and Reduction (2.2) to chains
  of morphisms

  \vspace{-1ex}
  {\footnotesize
  $$
   \coprod_{i,j}\,V_{i,j}^{\prime}\;
     \stackrel{\coprod_i\,f_i}{\longrightarrow}\;
     \coprod_i\, V_i\;\longrightarrow\; U
      \hspace{1em}\mbox{\normalsize and}\hspace{1em}
   \coprod_{i,j}\,V_{i,j}^{\prime}\;\longrightarrow\; U^{\prime}\;
     \stackrel{f}{\longrightarrow}\;U\,,
  $$
  {\normalsize one}} 
  concludes that $(a)$ and $(b)$ hold for $f$.
 If $f$ is fppf, then $V_{i,j}:=f(V_{i,j}^{\prime})$ are open in $U$.
  For each $i$, $\{V_{i,j}\}_j$ is a Zariski open cover of $V_i$ and
  each morphism $V_{i,j}^{\prime}\rightarrow V_{i,j}$ is fpqc.
  Applying Reduction (2.1) and Reduction to $(2.2)$ to chains
  of morphisms

  \vspace{-1ex}
  {\footnotesize
  $$
   \coprod_j\,V_{i,j}^{\prime}\;
     \longrightarrow\;
     \coprod_j\, V_{i,j}\;\longrightarrow\; V_i\,,
      \hspace{2em}
   \coprod_{i,j}\,V_{i,j}^{\prime}\;
     \longrightarrow\;
     \coprod_i\, V_i\;\longrightarrow\; U\,,
      \hspace{1em}\mbox{\normalsize and}\hspace{1em}
   \coprod_{i,j}\,V_{i,j}^{\prime}\;\longrightarrow\; U^{\prime}\;
     \stackrel{f}{\longrightarrow}\;U\,,
  $$
  {\normalsize one}} 
  concludes that $(a)$ and $(b)$ hold for $f$.
 This concludes the sketch.

\noindent\hspace{14cm}$\Box$

\bigskip

Recall that $\QM_g$ is the category fibered in groupoids of
 quasistable curves over $(\Sch/S_0)$.

\bigskip

\noindent
{\bf Fact 2.1.8 [$\QM_g$ Artin].} ([Be1].) {\it
 $\QM_g$ is an Artin stack for all $g\ge 0$.
} 

\bigskip

\noindent
{\it Explanation.}
 This follows from [Be1$\,$: Preliminaries on prestable curves]
  since $\QM_g$ is an open substack of the Artin stack of prestable
  curves.

\noindent\hspace{14cm}$\Box$

\bigskip

\noindent
{\bf Corollary 2.1.9 [line bundle with a section].} {\it
 The category fibered in groupoids of line bundles with a section
  on quasistable curves over $S_0$-schemes is a stack.
} 

\bigskip

\noindent
{\it Proof.}
 Since $\QM_g$ is an Artin stack, a descent datum for quasistable
  curves descends effectively.
 If $\underline{f}:S_1\rightarrow S_2$ is a fppf morphism and
  ${\cal C}_2$ is a quasistable curve over $S_2$, then
  ${\cal C}_1 := \underline{f}^{\ast}{\cal C}_2
               = S_1\times_{S_2}{\cal C}_2$
  is a quasistable curve over $S_1$ and
  the morphism $f:{\cal C}_1\rightarrow {\cal C}_2$ from
  fibered product is also fppf.
 We shall apply Fact 2.1.7 and its proof to
  $f:{\cal C}_1\rightarrow{\cal C}_2$.

 Given a descent data $({\cal E}^{\prime},\tau)$ on ${\cal C}_1$
  with ${\cal E}^{\prime}$ invertible, let ${\cal E}$ be its descent
  on ${\cal C}_2$. Since $f^{\ast}{\cal E}\simeq {\cal E}^{\prime}$
  and $f$ is fppf, ${\cal E}$ must be invertible as well.
  $$
   \begin{array}{cccl}
     {\cal C}_1 & \stackrel{f}{\longrightarrow}  & {\cal C}_2  & \\
     \hspace{2ex}\downarrow\mbox{\scriptsize $\pi_1$}
       & & \hspace{2ex}\downarrow\mbox{\scriptsize $\pi_2$}    &  \\
     S_1 & \stackrel{\underline{f}}{\longrightarrow} & S_2
   \end{array}
  $$
 Thus Fact 2.1.7 remains true if quasi-coherent sheaves are replaced
  by invertible sheaves
  (cf.\ [Mu1$\,$: Theorem 90 (Hilbert-Grothendieck)]).

 In the construction of the descent of quasi-coherent sheaves and
  their morphisms for the affine case in the proof of Fact 2.1.7,
  one observes that if a global section $s^{\prime}$ is added to
  a descent datum $({\cal E}^{\prime},\tau)$ that satisfies the gluing
  condition given by $\tau$, then it must lies in $\Ker\Theta$ and
  hence descends to a global section $s$ in ${\cal E}$.
 Furthermore, the statements in Item $(1.b^{\ast})$ in that proof
  imply that such $s$ is unique. I.e.\ the descent remains
  effective with a section added to the data.
 Consequently, Statement (a) and Statement (b) in Fact 2.1.7
  hold for descent data of invertible sheaves with a section
  when $f$ is a faithfully flat morphism between affine schemes.
 Now observe that in the remaining part of the proof of Fact 2.1.7,
  precisely two things are used repeatedly$\,$:
  \begin{itemize}
   \item [(i)] \parbox{13ex}{\it going-down$\,:$} \
     \parbox[t]{26em}{effective descent by a morphism, for which
                Statement (a) and Statement (b) are known to hold,}

   \item [(ii)] \parbox{13ex}{\it going-up$\,:$} \
     \parbox[t]{26em}{pulling back a descent datum.}
  \end{itemize}
 Whenever a going-down is employed, the existence and uniqueness of
  descent global section are known to hold by earlier reductions from
  the affine case
 while a going-up takes any part of descent data to the unique
  corresponding part of descent data automatically.
 Thus the whole proof of Fact 2.1.7 goes through without change.
 This proves the lemma.

\noindent\hspace{14cm}$\Box$

\bigskip

\noindent
{\it Proof of Definition/Lemma 2.1.5.}
 Continuing the notations from previous discussions.
 Since the nonvanishing condition of weak $\Delta$-collections
  is an open condition, we only need to show that the data without
  this condition give a stack.

 \medskip

 \noindent
 $(a)$ {\it Descent of weak $\Delta$-collections.}
 For each $\rho\in\Delta(1)$ the existence and uniqueness of descent
  of the descent data
  $(L_{\rho}^{\prime},u_{\rho}^{\prime}\,;\tau_{\rho})$ on
  ${\cal C}_1$ to $(L_{\rho},u_{\rho}\,;\sigma_{\rho})$ on
  ${\cal C}_2$ follow from Corollary 2.1.9.
 Similarly for morphisms between two such descent data.
 Each
  $(\otimes_{\rho}{L^{\prime}_{\rho}}^{\otimes\langle m,n_{\rho}\rangle},
     \otimes_{\rho}{\tau_{\rho}}^{\otimes\langle m,n_{\rho}\rangle})$,
  as well as $({\cal O}_{{\cal C}_1},\Id)$, is a descent datum of line
  bundles on ${\cal C}_1$.
 Their descent on ${\cal C}_2$ are given by
  $(\otimes_{\rho}{L_{\rho}}^{\otimes\langle m,n_{\rho}\rangle},
    \otimes_{\rho}\,{\sigma_{\rho}}^{\otimes\langle m,n_{\rho}\rangle})$
  and $({\cal O}_{{\cal C}_2},\Id)$ respectively.
 Thus
  $c_m^{\prime}:
   \otimes_{\rho}{L^{\prime}_{\rho}}^{\otimes\langle m,n_{\rho}\rangle}
   \rightarrow {\cal O}_{{\cal C}_1}$,
  as an isomorphism between two descent data of line bundles,
  descends to a unique isomorphism
  $c_m:\otimes_{\rho}{L_{\rho}}^{\otimes\langle m,n_{\rho}\rangle}
       \rightarrow {\cal O}_{{\cal C}_2}$.
 This shows that a descent datum of weak $\Delta$-collection on
  quasistable curves descends effectively.

 \medskip

 \noindent
 $(b)$ {\it Descent of morphisms.}
 Since a descent datum of isomorphisms from a weak
  $\Delta$-collection to another is really that for line bundles.
 It descends uniquely.

 \medskip
 This concludes the proof.

\noindent\hspace{14cm}$\Box$

\bigskip

\subsection{$\AM_g(X)$ is an Artin stack.}

\noindent
{\bf Proposition 2.2.1 [$\AM_g(X)$ Artin].} {\it
 The A-twisted moduli stack $\AM_g(X)$ for genus $g$ quasistable
  curves is an Artin stack.
} 

\bigskip

\noindent
{\it Proof.}
We check the properties that the diagonal morphism needs to satisfy
 and construct an atlas for $\AM_g(X)$ via a relative construction.

\medskip

\noindent
$(a)$ {\it Representability, quasi-compactness, and separatedness
           of the diagonal morphism.}

\medskip

\noindent
These properties are reflected in the corresponding properties of
 the $\underline{\Isom}$-functor (cf.\ [G\'{o}$\,$: Sec.\ 2.2]),
 which we will now check.

\bigskip

\noindent $(a.1)$
{\it Representability of $\underline{\Isom}_U({\cal F}_1,{\cal F}_2)$.}
Let
 $U\in(\Sch/S_0)$ and
 ${\cal F}_1$, ${\cal F}_2\in \AM_g(X)(U)$ with the underlying
  quasistable curves over $U$ denoted by
  $\pi_1:{\cal C}_1\rightarrow U$, $\pi_2:{\cal C}_2\rightarrow U$
  respectively.
Let $\underline{\Isom}_U({\cal F}_1,{\cal F}_2)$ be the functor on
 $(\Sch/U)$ associating to each $U$-scheme
 $f:U^{\prime}\rightarrow U$ the set of $U^{\prime}$-isomorphisms
 from $f^{\ast}{\cal F}_1$ to $f^{\ast}{\cal F}_2$.
By passing through a standard limit, one may assume that the base
 schemes $U$ and $U^{\prime}$ are Noetherian and of finite type
 over $S_0$ (cf.\ [L-MB$\,$: Th\'{e}or\`{e}me (4.6.2.1), proof]).
There is a natural morphism of functors
 $\underline{\Isom}_U({\cal F}_1,{\cal F}_2)
   \rightarrow\underline{\Isom}_U({\cal C}_1,{\cal C}_2)$,
 whose fiber over an isomorphism
 $\varphi:f^{\ast}{\cal C}_1\rightarrow f^{\ast}{\cal C}_2$
 over $f:U^{\prime}\rightarrow U$ is the set of isomorphisms
 from $f^{\ast}{\cal F}_1$ to $\varphi^{\ast}f^{\ast}{\cal F}_2$
 on $f^{\ast}{\cal C}_1$.
From [Gro] and Fact 2.1.8  
 the functor $\underline{\Isom}_U({\cal C}_1,{\cal C}_2)$
 is represented by a scheme $\Isombf_U({\cal C}_1,{\cal C}_2)$
 quasi-compact and separated over $U$.
Let $h_0:\Isombf_U({\cal C}_1,{\cal C}_2)\rightarrow U$ be the
 natural morphism, then there is a canonical isomorphism
 $\Phi:h_0^{\ast}{\cal C}_1
      \stackrel{\sim}{\rightarrow}h_0^{\ast}{\cal C}_2$.
Denote $h_0^{\ast}{\cal C}_1$ over
 $Y_0:=\Isombf_U({\cal C}_1,{\cal C}_2)$ by $\widetilde{{\cal C}_1}$.
Consider the functor
 $\underline{\Isom}_{\widetilde{{\cal C}_1}/Y_0}(
         h_0^{\ast}{\cal F}_1,\Phi^{\ast}h_0^{\ast}{\cal F}_2)$.
From [Gro$\,$: Sec.\ 4]
 (also [L-MB$\,$: Th\'{e}or\`{e}me (4.6.2.1), proof]),
 if one just focuses on the part of collections of line bundles
 $\{L_{\rho}\}_{\rho}$ in weak $\Delta$-collections,
 then the $\underline{\Isom}$-functor is represented by an open
 subscheme $Y_2$ of a scheme $Y_1$ that is affine and of finite type
 over $Y_0$.
The additional data$\,$: sections $u_{\rho}$ and trivialization
 isomorphisms $c_m$ of
 $\otimes_{\rho}L_{\rho}^{\otimes\langle m,n_{\rho}\rangle}$,
 specify a locally closed subscheme $Y_3$ of $Y_2$.
Thus $\underline{\Isom}_{\widetilde{{\cal C}_1}/Y_0}(
         h_0^{\ast}{\cal F}_1,\Phi^{\ast}h_0^{\ast}{\cal F}_2)$
 is represented by $Y_3$ over $Y_0$.
In summary,
$$
 \begin{array}{ccccccccc}
  Y_3\;& \stackrel{h_3}{\longrightarrow}  & Y_2\;
       & \stackrel{h_2}{\longrightarrow}  & Y_1\;
       & \stackrel{h_1}{\longrightarrow}
       & Y_0=\Isombf_U({\cal C}_1,{\cal C}_2)\;
       & \stackrel{h_0}{\longrightarrow}  & U\,.                 \\
     & \mbox{\parbox{3em}{\tiny
                locally\newline closed \newline immersion}}     &
       & \mbox{\parbox{3em}{\tiny open \newline immersion}}     &
       & \mbox{\parbox{2.4em}{\tiny
               affine, \newline of finite \newline type}}       &
       & \mbox{\parbox{3em}{\tiny quasi-compact,\newline separated}}
       & \mbox{\tiny Noetherian}
 \end{array}
$$

Let $p:Y_3\rightarrow Y_0$ be the composition $h_1\circ h_2\circ h_3$,
then there is a canonical isomorphism
 $\Psi:p^{\ast}{\cal F}_1\rightarrow p^{\ast}\Phi^{\ast}{\cal F}_2$
 over $p^{\ast}h_0^{\ast}{\cal C}_1/Y_0$.

\bigskip

\noindent
{\bf Claim.} {\it
 The functor $\underline{\Isom}_U({\cal F}_1,{\cal F}_2)$
  is represented by $Y_3$.
} 

\bigskip

\noindent
{\it Proof.}
 One checks the functorial properties for a scheme that represents
  an $\underline{\Isom}$-functor.
 Let $g:U^{\prime}\rightarrow U$ be a $U$-scheme and
  $\widetilde{\gamma}:
    g^{\ast}{\cal F}_1\rightarrow g^{\ast}{\cal F}_2$
  be an isomorphism.
 Let $\gamma:g^{\ast}{\cal C}_1\rightarrow g^{\ast}{\cal C}_2$ be
  the underlying isomorphism of quasistable curves over $U^{\prime}$.
  Then one may rewrite $\widetilde{\gamma}$ as an isomorphism
  $\widetilde{\gamma}:g^{\ast}{\cal F}_1
               \rightarrow \gamma^{\ast}g^{\ast}{\cal F}_2$
  on $g^{\ast}{\cal C}_1$.
 Since $Y_0$ represents the functor
  $\underline{\Isom}_U({\cal C}_1,{\cal C}_2)$, there is
  a unique $U$-morphism $h:U^{\prime}\rightarrow Y_0$ such that
  $(\gamma:g^{\ast}{\cal C}_1\rightarrow g^{\ast}{\cal C}_2)$
  is the pullback of the canonical isomorphism
  $(\Phi:h_0^{\ast}{\cal C}_1\rightarrow h_0^{\ast}{\cal C}_2)$
  over $Y_0$.
 Since $h_0\circ h=g$, the data
  $(\widetilde{\gamma}:g^{\ast}{\cal F}_1
                  \rightarrow \gamma^{\ast}g^{\ast}{\cal F}_2)$
  is the same as
  $\widetilde{\gamma}:h^{\ast}h_0^{\ast}{\cal F}_1
           \rightarrow h^{\ast}\Phi^{\ast}h_0^{\ast}{\cal F}_2$.
 Since $Y_3$ is the scheme that represents the functor
   $\underline{\Isom}_{\widetilde{{\cal C}_1}/Y_0}(
        h_0^{\ast}{\cal F}_1,\Phi^{\ast}h_0^{\ast}{\cal F}_2)$,
  there is a unique $Y_0$-morphism
  $\widetilde{h}:U^{\prime}\rightarrow Y_3$ (i.e.\ a lifting of $h$)
  such that $\widetilde{\gamma}$ is the pullback of $\Psi$ via
  $\widetilde{h}$.
 This shows that
  $\underline{\Isom}_U({\cal F}_1,{\cal F}_2)=\Hom(\,\cdot\,,Y_3)$
  and hence $Y_3$ represents
  $\underline{\Isom}_U({\cal F}_1,{\cal F}_2)$.

\noindent\hspace{14cm $\Box$}

\bigskip

\noindent $(a.2)$ {\it Separatedness.}
Recall the schemes $Y_i$ and the morphisms $h_i$ between them
 from the discussion in Part $(a.1)$.
From Hartshorne [Ha]$\,$:
 $(i)$   affine morphisms (cf.\ $h_1$) are separated
         [Ha$\,$: II.\ Exercise 5.17(b)],
 $(ii)$  open and closed - and hence locally closed - immersions
         (cf.\ $h_2$, $h_3$) are separated
         [Ha$\,$: II.\ Corollary 4.6(a)], and
 $(iii)$ composition of separated morphisms is separated,
it follows that $Y_3\rightarrow U$ is separated since $h_0$
is separated.

\bigskip

\noindent $(a.3)$ {\it Quasi-compactness.}
Recall again from Hartshorne [Ha$\,$: II.\ Exercise 3.3(a)] that
 a morphism of schemes is of finite type if and only if
 it is locally of finite type and quasi-compact.
Since all the morphisms $h_i$ are of finite type, so do their
 composition $Y_3\rightarrow U$, which then must be quasi-compact.

\medskip

These together justify that the diagonal morphism is representable,
 quasi-compact, and separated.

\bigskip

\noindent
$(b)$ {\it Construction of an atlas.}

\medskip

\noindent
We follow a relative construction, which is completed in three steps.

\medskip

\noindent $(b.1)$ {\it Atlas $U_0$ for $\QM_g$.}
It follows from the discussion in [Be1] that an atlas for $\QM_g$
 can be chosen as follows.
Observe that for $C$ a quasistable curve of genus $g\ge 0$,
 the number of unstable components of $C$ is bounded strictly by
 $1$, for $g=0,\;1$ and $3g-3$, for $g\ge 2$,
 and the number of marked points on $C$ needed to stabilize all
 these components is bounded by $n_0=3$ for $g=0$, $1$ for $g=1$,
 and $3g-3$ for $g\ge 2$.
Let $\overline{\cal M}^q_{g,n_0}$ be the open substack of the
 Deligne-Mumford stack $\overline{\cal M}_{g,n_0}$ that consists
 of stable curves of genus $g$ with $n_0$ marked points such that
 when these marked points are forgotten, the underlying curves are
 quasistable.
Then $\overline{\cal M}^q_{g,n_0}$ is a Deligne-Mumford stack
 and the morphism $F:\overline{\cal M}^q_{g,n_0}\rightarrow\QM_g$
 induced by forgetting the marked points is representable, smooth,
 and surjective. Indeed if $V\in(\Sch/S_0)$ and
  $(\pi:{\cal C}\rightarrow V)\in \QM_g(V)$,
  then the fibered-product morphism
  $\overline{\cal M}^q_{g,n_0}\times_{F,
        {\cal Q}{\cal M}_g,\pi}\,V \rightarrow V$
 is the morphism of $S_0$-schemes
 $$
  (\;\underbrace{{\cal C}\times_V\,\cdots\,\times_V {\cal C}}_{
             \mbox{\scriptsize $n_0$-many times}}\;)^{(0)}\;
  \longrightarrow\; V\,,
 $$
 where $({\cal C}\times_V\,\cdots\,\times_V {\cal C})^{(0)}$
  is the fibered product ${\cal C}\times_V\,\cdots\,\times_V {\cal C}$
  with all the diagonals of the fibered product $\times_V$ and
  the locus of nodes of fibers of ${\cal C}\rightarrow V$ removed.
This re-justifies that $\QM_g$ is an Artin stack.
Recall the Hilbert scheme construction in [$\,$F-P$\,$: Sec.\ 2$\,$]
 (with ${\Bbb P}^r$ therein set to ${\Bbb P}^0=\{\pt\}$)
 that realizes $\overline{\cal M}_{g,n_0}$ as a quotient stack
 of a quasi-projective variety $U$ acted on by an algebraic group.
(Caution that a ``quasistable curve" in [F-P] is a ``prestable curve"
 in [Be1] and the current article.)
There is an open subset $U_0$ in $U$ whose geometric points corresponds
 to quasistable curves.
This $U_0$ is then an atlas for $\QM_g$.
Since it comes from a Hilbert scheme construction, the associated
 flat family of quasistable curves
 $\pi:{\cal C}_0\rightarrow U_0$ is projective.
We shall fix a relative very ample line bundle on ${\cal C}_0/U_0$.

\bigskip

\noindent $(b.2)$ {\it Atlas $U_1$ for $\WD^X_{{\cal C}_0/U_0}$.}
Consider now the stack $\WD^X_{{\cal C}_0/U_0}$ over
 $(\Sch/S_0)$ of weak $\Delta$-collections on ${\cal C}_0/U_0$.
Define an atlas $V$ for a stack ${\cal S}$ in the same way as that
 for an Artin stack, namely a morphism $V\rightarrow{\cal S}$ that
 is representable, smooth, and surjective.
Then an atlas for $\WD^X_{{\cal C}_0/U_0}$ can be constructed by
 a sequence of relative constructions given in the following steps.
\begin{itemize}
 \item [] \hspace{-1.6em}\parbox{3.4em}{($b.2.1$)}
   The Quot-scheme construction for an atlas $V_1$ for the Artin stack
    $\Bun_1({\cal C}_0/U_0)$ of line bundles on ${\cal C}_0$ flat
    over $U_0$
    (cf.\ [G\'{o}mez$\,$: Sec.\ 2.3, Example 2.24] and
          [L-MB$\,$: Example (4.6.2)]).
   $V_1$ can be decomposed into a disjoint union of components labelled
    by the degree of the line bundles on fiber of
    ${\cal C}_0\rightarrow U_0$ since the degree determines
    the Hilbert polynomial of line bundles on curves of fixed genus.
   By construction there are
    a natural quasi-projective morphism $V_1\rightarrow U_0$ and
    a tautological line bundle $\widetilde{\cal L}$ on
     $V_1\times_{U_0}{\cal C}_0$ over $V_1$.

 \item [] \hspace{-1.6em}\parbox{3.4em}{($b.2.2$)}
  An atlas $V_2$ for the stack $\Bun_{1,s}({\cal C}_0/U_0)$
   of line bundles with a section on ${\cal C}_0/U_0$ is given by
   the scheme
   $V_2:=\Hombf_{(V_1\times_{U_0}{\cal C}_0)/V_1}(
         {\cal O}_{V_1\times_{U_0}{\cal C}_0}\,,\,\widetilde{\cal L})$.
   It is affine and surjective over $V_1$,
    ([Gro$\,$: Sec.\ 4] and [L-MB$\,$: Sec.\ (4.6.2)]).
  The fibered product
   $$
    V_3\;:=\;\underbrace{V_2\times_{U_0}\,\cdots\,\times_{U_0}V_2}_{
                   \mbox{\scriptsize $|\Delta(1)|$-many times}}
   $$
  gives an atlas for the stack of $|\Delta(1)|$-tuple of lines
   with a section on the same quasistable curves.

 \item [] \hspace{-1.6em}\parbox{3.4em}{($b.2.3$)}
  The tensor product conditions
   $\otimes_{\rho}\,L_{\rho}^{\otimes\langle m,n_{\rho}\rangle}
    \simeq {\cal O}_C$
   in a weak $\Delta$-collection on a quasistable curve $C$ are
   locally closed conditions.
  Together they determine a locally closed subscheme $V_4$ in $V_3$.
  Fix a basis for the $M$-lattice, then over
   $V_4\times_{U_0}{\cal C}_0$ there is a $\rank M$-tuple of
    line bundles $(\widetilde{\cal L}_m)_m$ defined by
    $(\otimes_{\rho}\,L_{\rho}^{\otimes\langle m,n_{\rho}\rangle}
      )_m$, where $m$ runs over the fixed basis.
  To add in the data of the choices
   of trivialization
   $c_m:\otimes_{\rho}\,L_{\rho}^{\otimes\langle m,n_{\rho}\rangle}
        \simeq {\cal O}_C$,
   one takes the scheme
    $V_5:=\oplus_{m\in {\rm basis}}\,
      \Isombf_{(\,V_4\times_{U_0}{\cal C}_0)/V_4}(
         \widetilde{\cal L}_m\:,\:
                       {\cal O}_{V_4\times_{U_0}{\cal C}_0}\,)$,
    which is an affine bundle over $V_4$ with fiber the abelian group
    $\prod^{{\it rank}\,M }\Spec k[t,t^{-1}] $.

 \item [] \hspace{-1.6em}\parbox{3.4em}{($b.2.4$)}
  Finally, let $U_1$ be the open subscheme in $V_5$ that corresponds
   to the nonvanishing condition of a weak $\Delta$-collection.
  Note that we start with $V_1$ that is smooth and surjective over
   the stack $\Bun_1({\cal C}_0/U_0)$.
   As we start to enlarge $\Bun_1({\cal C}_0/U_0)$ to tuples of line
   bundles, choice of sections, and so on or doing restriction by
   imposing open, closed or locally closed conditions, these extra data
   or restrictions do not have nontrivial automorphisms. Thus they do
   not influence the representability and the smoothness of the morphism
   of resulting $V_i$ to the related stack in the discussion.
  Also, by construction they are surjective.
  Thus $U_1$ is an atlas for the stack $\WD^X_{{\cal C_0}/U_0}$.
  By construction there is a natural morphism $U_1\rightarrow U_0$.
\end{itemize}

\bigskip

\noindent $(b.3)$ {\it $U_1$ as an atlas for $\AM_g(X)$.}
By construction there is a relative weak $\Delta$-collection
 ${\cal F}_1$ on the quasistable curve ${\cal C}_1/U_1$.
This gives a morphism $f_1:U_1\rightarrow \AM_g(X)$.
Let ${\cal F}$ be a relative weak $\Delta$-collection on
 a quasistable curve ${\cal C}_W$ over $W\in(\Sch/S_0)$.
 This specifies a morphism $f_W:W\rightarrow\AM_g(X)$.
By the functorial properties of $\Isombf$ and the morphism
 $U_1\rightarrow U_0$, one has the following natural morphisms
 $$
  \Isombf({\cal C}_1/U_1,{\cal C}_W/W)\;
   =\; U_1\times_{U_0}\Isombf({\cal C}_0/U_0,{\cal C}_W/W)\;
   \stackrel{\pi}{\longrightarrow}\;
       \Isombf({\cal C}_0/U_0,{\cal C}_W/W)
 $$
 and
 $$
  \begin{array}{c}
    \pi_1^{\ast}{\cal C}_1\simeq\pi^{\ast}\pi_{10}^{\ast}\,{\cal C}_0
      \hspace{1ex}\stackrel{\stackrel{\alpha\,:=\,\pi^{\ast}(\alpha_0)
                                      }{\sim}}{\longrightarrow}
       \hspace{1ex}
      \pi_2^{\ast}\,{\cal C}_W\simeq\pi^{\ast}\pi_{20}^{\ast}{\cal C}_W\\
    \setminus \hspace{2em} / \\
    \hspace{7em}\Isombf({\cal C}_1/U_1,{\cal C}_W/W)\\
    \mbox{\scriptsize $\pi_1$}\swarrow \hspace{1em}
      \searrow \mbox{\scriptsize $\pi_2$}\\
    U_1 \hspace{3em} W
  \end{array}
  \hspace{1em} \stackrel{\pi}{\longrightarrow}\hspace{-3em}
  \begin{array}{c}
    \pi_{10}^{\ast}\,{\cal C}_0 \hspace{1ex}
      \stackrel{\stackrel{\alpha_0}{\sim}}{\longrightarrow}
      \hspace{1ex} \pi_{20}^{\ast}\,{\cal C}_W \\
    \setminus \hspace{2em} / \\
    \hspace{7em}\Isombf({\cal C}_0/U_0,{\cal C}_W/W)\,. \\
    \mbox{\scriptsize $\pi_{10}$}\swarrow \hspace{1em}
      \searrow \mbox{\scriptsize $\pi_{20}$}\\
    U_0 \hspace{3em} W
  \end{array}
 $$
It follows that
 \begin{eqnarray*}
  \lefteqn{ U_1\times_{f_1,\,{\cal A}{\cal M}_g(X),\,f_W}W\;
   =\; \Isombf({\cal F}_{{\cal C}_1/U_1}, {\cal F}_{{\cal C}_W/W})} \\
   && =\;
       \Isombf_{\pi_1^{\ast}{\cal C}_1/
                    {\bf Isom}({\cal C}_1/U_1,\,{\cal C}_W/W)}\,
       (\pi_1^{\ast}{\cal F}_1,\alpha^{\ast}\pi_2^{\ast}\,{\cal F}_W) \\
   && =\;
       \Isombf\left(\,({\cal F}_1)_{{\cal C}_1/U_1/U_0}\,,\;
          (\alpha_0^{\ast}\pi_{20}^{\ast}
          {\cal F}_W)_{(\pi_{10}^{\ast}\,{\cal C}_0)/U_0}\,\right) \\
   && =\; U_1 \times_{{\cal WD}^X_{{\cal C}_0/U_0}}
              \Isombf({\cal C}_0/U_0,\, {\cal C}_W/W) \\
   && \stackrel{\mbox{\scriptsize smooth and surjective}
                }{\longrightarrow}\hspace{1em}
      \Isombf({\cal C}_0/U_0,\, {\cal C}_W/W)\;
      =\; U_0 \times_{{\cal QM}_g} W \\
   && \stackrel{\mbox{\scriptsize smooth and surjective}
                }{\longrightarrow}\hspace{1em} W\,.
 \end{eqnarray*}
This shows that $U_1$ is an atlas for $\AM_g(X)$ and
we conclude the proof.

\noindent\hspace{14cm}$\Box$

\bigskip

\noindent
{\it Remark 2.2.2.}
 The above type of relative construction can be found also in the
 study of relative GIT construction of universal moduli spaces,
 e.g., [Hu] and [Pa].

\bigskip

\bigskip

\section{The $\;g=0\;$ case.}

Since $\CP^1$ is rigid, the problem may be treated as in the study
 of bundles on a fixed variety.
$\AM_0(X)$ is then the stackification of the prestack
 $\pre\AM_0(X)$, whose fiber $\pre\AM_0(X)$ over
 $S\in(\Sch/S_0)$ is the groupoid
 $$
  \pre\AM_0(X)(S)\;
  =\;\{\,\mbox{weak $\Delta$-collections
                            on $S\times\CP^1$ over $S$}\,\}\,.
 $$

On the other hand one has the construction of Morrison-Plesser
 [M-P$\,$: Sec.\ 3.7], as is used in
 [L-L-Y$\,$: II, Sec.\ 2.4, Example 4].
In this section we shall discuss how Morrison-Plesser's construction
 is related to Cox's work and the stack $\AM_0(X)$ adapted
 from Sec.\ 2.
We shall assume that $X=X_{\Delta}$ is convex throughout this section.
In particular, this implies that every entry $d_{\rho}$ of
 a multi-degree $d=(d_{\rho})_{\rho}$ in the discussions are
 all nonnegative integers.

\bigskip

\begin{flushleft}
{\bf The small universal weak $\Delta$-collection
     \`{a} la Morrison-Plesser.}
\end{flushleft}
Fix a presentation$\,$:
 ${\Bbb P}^1=\Proj {\Bbb C}[z_0,z_1]\,$,
 ${\Bbb C}[z_0,z_1]=\oplus_{l\ge 0} M_l\,,$
 ${\cal O}_{{\scriptsizeBbb P}^1}(l)={\Bbb C}[z_0,z_1](l)^{\,\sim}$
 and
 $H^0({\Bbb P}^1,{\cal O}_{{\scriptsizeBbb P}^1})(l)=M_l$
  for $l\ge 0\,$, cf. [Ha1].
Then the graded ${\Bbb C}$-algebra structure
 $M_{l_1}\cdot M_{l_2}\rightarrow M_{l_1+l_2}$
 induces a set of canonical isomorphisms of sheaves
 ${\cal O}_{{\scriptsizeBbb P}^1}(l_1)
  \otimes{\cal O}_{{\scriptsizeBbb P}^1}(l_2)
  \rightarrow{\cal O}_{{\scriptsizeBbb P}^1}(l_1+l_2)$.
Since multiplication among $M_l$'s is associative with respect to
 these isomorphisms, one has also canonical isomorphisms
 ${\cal O}(l_1)\otimes\cdots\otimes{\cal O}(l_s)
  \rightarrow{\cal O}(l_1+\cdots+l_s)\,$.
This implies that if let
  $\Delta(1)=\{\,\rho_{1i}\,\}_i\coprod\{\,\rho_{2j}\,\}_{j}$
  such that $\sum_i m_{1i}d_{\rho_{1i}}=\sum_j m_{2j}d_{\rho_{2j}}$
   for some $m_{1i}$, $m_{2j}\,\ge 0$,
 then there is a canonical isomorphism
  $\bigotimes_i\,
    {\cal O}_{{\scriptsizeBbb P}^1}(d_{\rho_{1i}})^{\;\otimes m_{1i}}
   \simeq
   \bigotimes_j\,
   {\cal O}_{{\scriptsizeBbb P}^1}(d_{\rho_{2j}})^{\;\otimes m_{2j}}\,$.
 (Cf.\ [Cox1$\,$: Proposition 1.1].)

\bigskip

\noindent
{\bf Definition 3.1 [set of canonical isomorphisms].} {\rm
We shall call the above set of isomorphisms
 the set of {\it canonical isomorphisms} among tensor products of
 ${\cal O}_{{\scriptsizeBbb P}^1}(l)$'s
 {\it with respect to the fixed presentation}.
} 

\bigskip

Let
 $d=(d_{\rho})_{\rho\in\Delta(1)}$, with $d_{\rho}$ nonnegative
 integers, be a multi-degree and define
 $$
  Y_d\;:=\;\oplus_{\rho}\,H^0({\Bbb P}^1,{\cal O}(d_{\rho}))
            \hspace{1em}\mbox{with the above fixed presentation}\,.
 $$
Recall the abelian group $G$ and the quotient
 $X=({\Bbb C}^{\Delta(1)}-V(I))/\mbox{\raisebox{-.4ex}{$G$}}$
 from Explanation/Fact 2.1.2.
Then each element of $Y_d$ corresponds to a morphism
 ${\Bbb P}^1\rightarrow {\Bbb C}^{\Delta(1)}$ up to
 a ${\Bbb C}^{\times}$-action on ${\Bbb C}^{\Delta(1)}$ by
 $t\cdot (x_{\rho})_{\rho}=(t^{\,d_\rho}x_{\rho})_{\rho}$.
Define $F_d$ to be the subvariety of $Y_d$ that consists of elements
 whose corresponding map ${\Bbb P}^1\rightarrow{\Bbb C}^{\Delta(1)}$
 has image contained in $V(I)\,$.
Since $V(I)$ is a union of coordinate subspaces in ${\Bbb C}^{\Delta(1)}$
 and hence invariant under the above ${\Bbb C}^{\times}$-action,
 $F_d$ is well-defined.
The $G$-action on ${\Bbb C}^{\Delta(1)}$ induces a $G$-action on $Y_d$
 that leaves $F_d$ invariant. Thus, one can define the quotient
 $$
  W_d\;=\;{\cal M}_d\,:=\,(Y_d-F_d)/G\,.
 $$
Cf.\ Appendix, [M-P$\,$: Sec.\ 3.7];
     also [L-L-Y$\,$: II, Sec.\ 2.4, Example 4].

Let ${\cal F}:=(L_{\rho},u_{\rho},c_m)_{\rho,m}$
 be a weak $\Delta$-collection on ${\Bbb P}^1$ of multi-degree $d$.
Then $(c_m)_m$ determines isomorphisms
 $L_{\rho}\simeq{\cal O}_{{\scriptsizeBbb P}^1}(d_{\rho})$,
 compatible with the set of canonical isomorphisms,
 up to an ambiguity parameterized by $G$.
Thus ${\cal F}$ corresponds to a $G$-orbit $O_{\cal F}$ in $Y_d\,$.
The nonvanishing condition for ${\cal F}$ is that
 $\;\sum_{\sigma\in\Delta_{max}}\otimes_{\rho\not\subset\sigma}
    u_{\rho}^{\ast}\,:\,
   \oplus_{\sigma\in\Delta_{max}}\otimes_{\rho\not\subset\sigma}
    L_{\rho}^{-1}\, \rightarrow\, {\cal O}_{{\scriptsizeBbb P}^1}\;$
 is not a zero-morphism.
Since $V(I)$ is defined by the ideal
 $I=(\,\prod_{\rho\not\subset\sigma}x_{\rho}\,|\,
                                     \sigma\in\Delta_{\it max}\,)$
 and the divisor on ${\Bbb C}^{\Delta(1)}$ defined by $x_{\rho}$
 corresponds to the subscheme on ${\Bbb P}^1$ defined by $u_{\rho}$,
 the nonvanishing condition means precisely that
 $O_{\cal F}\subset Y_d-F_d\,$.

Regard the sections of ${\cal O}_{{\scriptsizeBbb P}^1}(d_{\rho})$
 as subschemes of the total space
 $\Spec\Sym^{\bullet}({\cal O}_{{\scriptsizeBbb P}^1}(d_{\rho})^{\vee})$
 of ${\cal O}_{{\scriptsizeBbb P}^1}(d_{\rho})$.
Then as in the case of Hilbert schemes one obtains the universal
 tuple of sections $(\widetilde{u}_{\rho})_{\rho}$ of the line bundles
 $(\,\widetilde{\cal O}(d_{\rho})\,)_{\rho}$ over
 $(Y_d-F_d)\times{\Bbb P}^1$ from the pullback of the projection map
 $(Y_d-F_d)\times{\Bbb P}^1\rightarrow {\Bbb P}^1$.
The set of canonical isomorphisms in Definition 3.1
 gives a canonical
 set of isomorphisms $(\widetilde{c}_m)_m$.
Since $Y_d-F_d$ corresponds to tuples $(u_{\rho})_{\rho}$ of sections
 that satisfy the nonvanishing condition,
 $(\widetilde{\cal O}(d_{\rho}),
   \widetilde{u}_{\rho}, \widetilde{c}_m)_{\rho,m}$
 is a weak $\Delta$-collection on $(Y_d-F_d)\times{\Bbb P}^1$
 over $Y_d-F_d$.

\bigskip

\noindent
{\bf Definition 3.2 [universal weak $\Delta$-collection].} {\rm
 $\widetilde{\cal F}_d
   :=(\,\widetilde{\cal O}(d_{\rho}),\widetilde{u}_{\rho},
        \widetilde{c}_m\,)_{\rho,m}$
 is called the {\it universal weak $\Delta$-collection} on
 $(Y_d-F_d)\times{\Bbb P}^1$ over $Y_d-F_d\,$.
} 

\bigskip

\noindent
{\it Remark 3.3 {\rm [}universal property{\rm ]}.}
 Since the line bundles $\widetilde{\cal O}(d_{\rho})$ in
  $\widetilde{\cal F}_d$ are fixed and $\widetilde{c}_m$
  is determined once a presentation of ${\Bbb P}^1$ is chosen,
  $\widetilde{\cal F}_d$ indeed comes from a restriction of
  the universal subscheme over a Hilbert scheme.
 It thus inherits a similar universal property as Hilbert schemes.

\bigskip

The $({\Bbb C}^{\times})^{|\Delta(1)|}$-action on $Y_d-F_d$
 lifts to a $({\Bbb C}^{\times})^{|\Delta(1)|}$-action on
 $(Y_d-F_d)\times{\Bbb P}^1$ by acting on ${\Bbb P}^1$ by
 the identity.
The latter then lifts to an action on each
 $\widetilde{\cal O}(d_{\rho})$ that leaves $\widetilde{u}$ invariant.
This is the unique lift that has this property.
Since $G$ is a subgroup of $({\Bbb C}^{\times})^{|\Delta(1)|}$,
 $G$ lifts to a unique action on
 $(\widetilde{\cal O}(d_{\rho}),\widetilde{u}_{\rho})_{\rho}$ as well.
By the very definition of $G$, this $G$-action commutes with
 $(\widetilde{c}_m)_m\,$.

\bigskip

\noindent
{\bf Definition 3.4 [canonical $G$-action].} {\rm
 The above $G$-action on
  $(\widetilde{\cal O}(d_{\rho}),\widetilde{u}_{\rho},
                               \widetilde{c}_m)_{\rho,m}$
 is called the {\it canonical lifting} of the $G$-action on $Y_d-F_d$.
} 

\bigskip

\begin{flushleft}
{\bf The (big) universal weak $\Delta$-collection \`{a} la Cox.}
\end{flushleft}
Continuing the notations in the previous theme.
Fix a basis of the $M$-lattice. Define
 $$
  \Xi_{\,d}\;
   =\; \oplus_m\Isom_{{\cal O}_{{\tinyBbb P}^1}}(\,
        \otimes_{\rho}{\cal O}_{{\scriptsizeBbb P}^1}
                  (d_{\rho})^{\otimes\langle m,n_{\rho}\rangle}\,,\,
         {\cal O}_{{\scriptsizeBbb P}^1}\,)
  \hspace{1em}\mbox{and}\hspace{1em}
  \widehat{Y}_d\;
   :=\; Y_d\oplus\Xi_{\,d}\,,
 $$
 where $m$ runs over the fixed basis of $M$.
Let $\kappa:\widehat{Y}_d\rightarrow \Xi_{\,d}$ be the natural
 projection.
Then, similar to the discussion in the previous theme, one has the
 $({\Bbb C}^{\times})^{\Delta(1)}
  :=\Hom_{\scriptsizeBbb Z}({\Bbb Z}^{\Delta(1)},
  {\Bbb C}^{\times})$-action
 on $\widehat{Y}_d$ induced from that on ${\Bbb C}^{\Delta(1)}$.
Recall $T_N$ from Definition/Fact 2.1.2,
 then $T_N$ acts on $\Xi_{\,d}$ freely and transitively and
 in such a way that $\kappa$ is
 $(({\Bbb C}^{\times})^{\Delta(1)},T_N)$-equivariant
 and that the $G$-subaction on $\widehat{Y}_d$ leaves each preimage
 of $\kappa$ invariant.
Since the nonvanishing condition in a weak $\Delta$-collection has
 nothing to do with the isomorphism data $(c_m)_m$, it specifies
 the open $({\Bbb C}^{\times})^{\Delta(1)}$-invariant subvariety
 $\widehat{Y}_d-\widehat{F}_d$,
 where $\widehat{F}_d:= F_d\times\Xi_{\,d}$.

Following the same construction as in the previous theme,
 one has a universal weak $\Delta$-collection
 $$
  \widetilde{\cal F}^{\rm\, big}_d\;
   :=\; \left(\,\widetilde{\cal O}(d_{\rho})^{\rm big}\,,\,
        \widetilde{u}_{\rho}^{\rm\,big}\,,\,
        \widetilde{c}_m^{\rm\,big}\,\right)_{\rho,m}
 $$
 on
 $((\widehat{Y_d}-\widehat{F}_d)\times{\Bbb P}^1)/
    (\widehat{Y}_d-\widehat{F}_d)$
 and the canonical lifting of
 the $({\Bbb C}^{\times})^{\Delta(1)}$-action on the total space of
 line bundles $(\widetilde{\cal O}(d_{\rho})^{\rm big})_{\rho}$.

\bigskip

\noindent
{\it Remark 3.5 {\rm [}Morrison-Plesser v.s.\ Cox$\,${\rm ]}.}
 From these very explicit constructions, one observes that a fixed
  presentation as in the previous theme selects a distinguished
  point $(c_m^{\rm can})_m$ in $\Xi$ and
  $Y_d-F_d=\kappa^{-1}((c_m^{\rm can})_m)\,$.
 The universal weak $\Delta$-collection $\widetilde{\cal F}_d$
  on $((Y_d-F_d)\times{\Bbb P}^1)/(Y_d-F_d)$
  is the restriction of $\widetilde{\cal F}^{\rm\, big}_d$ on
  $((\widehat{Y}_d-\widehat{F}_d)\times{\Bbb P}^1)/
                             (\widehat{Y}_d-\widehat{F}_d)$
  to $((Y_d-F_d)\times{\Bbb P}^1)/(Y_d-F_d)\,$.
 The $\Isombf$ construction in Sec.\ 2.2, adjusted for the fixed
  ${\Bbb P}^1$, gives $\widehat{Y}_d-\widehat{F}_d$.

\bigskip

\begin{flushleft}
{\bf Relation with $\AM_0(X)$.}
\end{flushleft}
The two quotient stacks
 $[(\widehat{Y}_d-\widehat{F}_d)/({\Bbb C}^{\times})^{\Delta(1)}]$ and
 $[(Y_d-F_d)/G]$ are isomorphic since the $T_N$-action on $\Xi_{\,d}$
 is transitive and free.
The following lemma relates this quotient stack with $\AM_0(X)$.

\bigskip

\noindent
{\bf Lemma 3.6 [$\AM_0(X)$].} {\it
 The Artin stack $\AM_0(X)$ is the quotient stack
  $\coprod_d\left[(Y_d-F_d)/G\right]$, for which
  $\coprod_d(Y_d-F_d)$ is an atlas and $\coprod_d\,W_d$ is
  the coarse moduli space.
 In particular, $\AM_0(X)$ is a smooth Artin stack.
} 

\bigskip

We check this at the prestack level. The statement then follows upon
 stackification.

\bigskip

\noindent
{\it Proof.} The proof is divided in two parts.

\medskip

\noindent $(a)$
{\it $\AM_0(X)$ as a quotient stack.}
Let $\pre\AM_0(X)=\coprod_d \pre\AM_0(X)_d$, where $d$ runs
 over all the admissible multi-degrees.
We shall construct morphisms of prestacks
 $$
  J^{(1)}_d:\pre\AM_0(X)_d\longrightarrow \pre[(Y_d-F_d)/G]
 $$
 and
 $$
  J^{(2)}_d:\pre[(Y_d-F_d)/G]\longrightarrow \pre\AM_0(X)_d
 $$
 so that $J^{(2)}_d\circ J^{(1)}_d$ and $J^{(1)}_d\circ J^{(2)}_d$
 induce auto-equivalences of related fiber groupoids.
 (In other words, $J^{(1)}_d$ is an isomorphism of prestacks with
  inverse given by $J^{(2)}_d$.)

\bigskip

\noindent
$(a.1)$ {\it Construction of $J^{(1)}_d$.}
Given a weak $\Delta$-collection
 ${\cal F}=(L_{\rho},u_{\rho},c_m)_{\rho,m}$ on $(S\times{\Bbb P}^1)/S$,
 each $L_{\rho}$ defines a principal ${\Bbb C}^{\times}$-bundle
 $L_{\rho}^{\times}$ over $S\times{\Bbb P}^1$ by deleting the
 zero-section of $L_{\rho}$.
The isomorphism of line bundles
 $c_m:\otimes_{\rho}L_{\rho}^{\otimes\langle m,n_{\rho}\rangle}
       \simeq {\cal O}_{S\times{\scriptsizeBbb P}^1}$
 induces an isomorphism of principal ${\Bbb C}^{\times}$-bundles over
 $S\times{\Bbb P}^1$ by the restriction
 $\otimes_{\rho}(L_{\rho}^{\times})^{\otimes\langle m,n_{\rho}\rangle}\;
   \rightarrow {\cal O}_{S\times{\scriptsizeBbb P}^1}^{\times}
               =({\Bbb G}_m)_{S\times{\scriptsizeBbb P}^1}$
 of $c_m$.
This then induces a morphism, still denoted by $c_m$,
 from the composition
 $$
  c_m\,:\,\oplus_{\rho}L_{\rho}^{\times}\;
   \longrightarrow\;
    \otimes_{\rho}(L_{\rho}^{\times})^{\otimes\langle m,n_{\rho}\rangle}\;
    \longrightarrow\:{\cal O}^{\times}_{S\times{\scriptsizeBbb P}^1}\,.
 $$
This gives rise to a principal $G$-bundle on $S\times{\Bbb P}^1$
 defined by the kernel (i.e.\ the preimage of the section
 $(1,\ldots,1)$) of the morphism $(c_m)_m$ over $S\times{\Bbb P}^1\,$:
 $$
  \Ker\left(\,(c_m)_m: \oplus_{\rho}\,L_{\rho}^{\times} \rightarrow
    ({\cal O}_{S\times{\scriptsizeBbb P}^1}^{\times})^{\,\oplus\,n}
    =(T_N)_{S\times{\scriptsizeBbb P}^1}\, \right)\,,
 $$
 where $m$ runs over elements in a fixed basis of $M$ and $n$
 is the rank of $M$.
A principal $G$-bundle over $S$, $p:P^G_S\rightarrow S$, is obtained
 by restricting the above principal $G$-bundle over $S\times{\Bbb P}^1$
 to a horizontal slice, e,g.\ $S\times\{0\}$.
(Note that any two such restrictions are isomorphic.
  The inverse of any such restriction of $L_{\rho}$ gives the line
  bundles on $S$ needed to twist $L_{\rho}$ so that the result is
  a pullback line bundle from that on ${\Bbb P}^1$.)
Consider the pullback weak $\Delta$-collection $p^{\ast}{\cal F}$
 on $(P^G_S\times{\Bbb P}^1)/P^G_S$.
The identity morphism of line bundles
 $(L_{\rho}|_{S\times\{0\}})_{\rho}
   \rightarrow(L_{\rho}|_{S\times\{0\}})_{\rho}$
 specifies a canonical trivialization
  $p^{\ast}(L_{\rho}|_{S\times\{0\}})\simeq {\cal O}_{P^G_S\times\{0\}}$
  on the horizontal slice $P^G_S\times\{0\}$ of $P^G_S\times{\Bbb P}^1$
  over $P^G_S$.
This implies that $p^{\ast}(L_{\rho})_{\rho}$ is isomorphic to
 the pullback of $({\cal O}_{{\scriptsizeBbb P}^1}(d_{\rho}))_{\rho}$
 by the projection map $P^G_S\times{\Bbb P}^1\rightarrow {\Bbb P}^1$.

Furthermore, a generalization of the following standard
 constructions$\,$:
 \begin{itemize}
   \item []
    Let $L^{\times}$ be the principal ${\Bbb C}^{\times}$-bundle
     on $S$ from deleting the zero-section of a line bundle $L$
     on $S$.
    The projection map $L^{\times}\rightarrow S$ pulls back
     $L$ to a line bundle $\widetilde{L}$ on the total space
     $\Tot(L^{\times})$ of $L^{\times}$.
    The natural inclusion map $L^0\hookrightarrow L$ gives rise to
     a nowhere-zero global section in $\widetilde{L}$ over
     $\Tot(L^{\times})$ and hence a canonical trivialization of
     $\widetilde{L}$.
 \end{itemize}
 to the tuple of line bundles $(L_{\rho}|_{S\times\{0\}})_{\rho}$,
 its associated principal $({\smallBbb C}^{\times})^n$-bundles
 and its sub $G$-bundle,
one deduces that the above trivialization over the slice
 $P^G_S\times\{0\}$ fixes an isomorphism
 $(p^{\ast}L_{\rho})_{\rho}
  \simeq ({\cal O}_{{\scriptsizeBbb P}^1}(d_{\rho}))_{\rho}$.

Pulling back now the sections $u_{\rho}$ of $L_{\rho}$, one thus
 obtains a $P^G_S$-family $p^{\ast}(L_{\rho},u_{\rho})_{\rho}$
 of line bundles on ${\Bbb P}^1$ with a section.
It follows from the universal property of $Y_d-F_d$ inherited
 from that of Hilbert schemes that there exists a unique morphism
  $\zeta_{\cal F}:P^G_S\rightarrow(Y_d-F_d)$
 with $p^{\ast}(L_{\rho},u_{\rho})_{\rho}
    =\zeta_{\cal F}^{\ast}\widetilde{\cal F}_d\,$.
By construction $\zeta_{\cal F}$ is $G$-equivariant and
 $p^{\ast}(c_m)_m=\zeta_{\cal F}^{\ast}\widetilde{c}_m$.
The correspondence ${\cal F}\rightarrow \zeta_{\cal F}$
 gives a morphism of prestacks
 $$
  J^{(1)}_d:\pre\AM_0(X)_d\rightarrow \pre[(Y_d-F_d)/G]\,.
 $$

\bigskip

\noindent
$(a.2)$ {\it Construction of $J^{(2)}_d$.}
To construct $J^{(2)}_d:\pre[(Y_d-F_d)/G]\rightarrow\pre\AM_0(X)$,
 observe that if $P^G_S$ is a principal $G$-bundle on $S$ with
 a $G$-equivariant morphism $\zeta:P^G_S\rightarrow (Y_d-F_d)$,
 (i.e.\ an object in the groupoid $\pre[(Y_d-F_d)/G](S)$)
 then $\zeta^{\ast}\widetilde{\cal F}_d$ is a weak
 $\Delta$-collection on $(P^G_S\times{\Bbb P}^1)/P^G_S$.
Since there is a canonical $G$-action on $\widetilde{\cal F}_d$,
 the $G$-action on $P^G_S$ also lifts canonically to
 $\zeta^{\ast}\widetilde{\cal F}_d$. The quotient by this action gives
 then a weak $\Delta$-collection ${\cal F}$ on $(S\times{\Bbb P}^1)/S$.
 i.e.\ an object in $\pre\AM_0(X)_d(S)$.
This gives a morphism
 $J^{(2)}_d:\pre[(Y_d-F_d)/G]\rightarrow\pre\AM_0(X)_d$.

\bigskip

\noindent
$(a.3)$ {\it Isomorphisms of stacks.}
It remains to show that $J^{(1)}_d$ (or $J^{(2)}_d$) is
 an isomorphism of stacks.
This means that $J^{(2)}\circ J^{(1)}$ sends a weak $\Delta$-collection
 on $(S\times{\Bbb P}^1)/S$ to an isomorphic weak $\Delta$-collection
 on $(S\times{\Bbb P}^1)/S$, which follows from the very explicit
 construction of $J^{(1)}_d$ and $J^{(2)}_d$.
Similarly for $J^{(1)}_d\circ J^{(2)}_d$.

\bigskip

\noindent $(b)$
{\it $W_d$ as the coarse moduli space.}

\medskip

\noindent
$(b.1)$ {\it Construction
             of a morphism $\AM_0(X)_d\rightarrow W_d$.}
A morphism $\pre\AM_0(X)_d\rightarrow W_d$ is already
 given/hidden in [L-L-Y$\,$: II. Sec.\ 2.5, Lemma 2.7, proof]
 as follows.
Let ${\cal O}(d_{\rho})$ be the pullback of
 ${\cal O}_{{\scriptsizeBbb P}^1}(d_{\rho})$ to $S\times{\Bbb P}^1$
 via the projection map.
Given a weak $\Delta$-collection
 ${\cal F}=(L_{\rho},u_{\rho}, c_m)_{\rho,m}$ on $S\times{\Bbb P}^1$
 over $S$ of multi-degree $d$,
 the data $(L_{\rho}, u_{\rho})_{\rho}$ determines non-uniquely a
 $(u_{\rho}^{\prime})_{\rho}
   \in \oplus_{\rho}H^0(S\times{\Bbb P}^1,{\cal O}(d_{\rho}))$
 by looking at the zero-divisor/locus of $u_{\rho}$
 on $S\times{\Bbb P}^1$.
The ambiguities are parameterized by
 $(t_{\rho})_{\rho}\in ({\Bbb C}^{\times})^{|\Delta(1)|}$
 that satisfy $\prod_{\rho} t_{\rho}^{\langle m,n_{\rho}\rangle}=1$
 for all $m\in M$.
Thus $(u_{\rho}^{\prime})_{\rho}$, though nonunique, determines
 a unique $S$-family of $G$-orbits on
 $\oplus_{\rho}H^0({\Bbb P}^1,{\cal O}(d_{\rho}))$.
The nonvanishing condition on $(u_{\rho}^{\prime})_{\rho}$
 inherits from that of ${\cal F}$.
Thus one obtains a morphism $S\rightarrow W_d$.
Another construction can be obtained from the discussion of Part (a)
 as follows.
 ${\cal F}$ determined a principal $G$-bundle $P^G_S$ on $S$
 with a unique $G$-equivariant morphism $P^G_S\rightarrow Y_d-F_d$.
 Taking quotient by $G$ on both sides, one then obtains a morphism
 $S\rightarrow W_d$ determined by ${\cal F}$.
Either way, one obtains a morphism
 $$
  \pre\,\phi\;:\; \pre\AM_0(X)_d\; \longrightarrow\; W_d
 $$
and hence a morphism
 $$
  \phi\;:\; \AM_0(X)_d\; \longrightarrow\; W_d
 $$
since for any $S\in(\Sch/S_0)$, $\pre\,\phi(S)$ depends only on the
 isomorphism class of the weak $\Delta$-collection on each fiber
 ${\Bbb P}^1$ of $S\times{\Bbb P}^1$ over $S$.

\bigskip

\noindent
$(b.2)$ {\it The coarse moduli space conditions.}
From the definition of $W_d$ and points of a stack,
 $|\phi|(k):|\AM_0(X)_d|(k)\rightarrow W_d(k)$
 is bijective for all algebraically closed field $k$.

To see that $W_d$ corepresents $\AM_0(X)_d$, observe
 that there is a distinguished weak $\Delta$-collection on
 $(W_d\times{\Bbb P}^1)/W_d$ constructed as follows.
Consider the diagonal action of $G$ on $(Y_d-F_d)\times(Y_d-F_d)$.
 The diagonal $\Delta_{(Y_d-F_d)}$ of $(Y_d-F_d)\times(Y_d-F_d)$
 is invariant under this $G$-action.
 The quotient gives a fibration of
 $\left((Y_d-F_d)\times(Y_d-F_d)\right)/G\rightarrow W_d$ with generic
 fiber $Y_d-F_d$. The diagonal $\Delta_{(Y_d-F_d)}$ descends to
 a section of this fibration, which corresponds to a weak
 $\Delta$-collection $\widetilde{\cal L}$ on $W_d\times{\Bbb P}^1$.
By construction $\phi(\widetilde{\cal L})=\Id_{\,W_d}$.
Suppose that $W_d^{\prime}$ is another scheme with a morphism
 $\phi^{\prime}:\AM_0(X)_d\rightarrow W_d^{\prime}$.
Define $\eta:\Hom(\,-\,,W_d)\rightarrow\Hom(\,-\,,W_d^{\prime})$
 by the composition
 $$
  (f:S\rightarrow W_d)\;
   \longmapsto\;
    f^{\ast}\widetilde{\cal L}\in\AM_0(X)_d(S)\;
   \longmapsto\;
    (\phi^{\prime}(f^{\ast}\widetilde{\cal L}):
                              S\rightarrow W_d^{\prime})\,.
 $$
Now given ${\cal L}\in\AM_0(X)_d(S)$, let
 $f=\phi({\cal L})\in\Hom(S,W_d)$ and
  $f^{\prime}=\phi^{\prime}({\cal L})\in \Hom(S,W_d^{\prime})$.
Then
 $$
  \eta(f)\;=\; \eta(\phi(f^{\ast}\widetilde{\cal L}))\;
   =\; \phi^{\prime}(f^{\ast}\widetilde{\cal L})\;
   =\; \phi^{\prime}({\cal L})\,,
 $$
 where we have used
 the observation that ${\cal L}$ and $f^{\ast}\widetilde{\cal L}$
 are fiberwise isomorphic weak $\Delta$-collections on
 $(S\times{\Bbb P}^1)/S$
 and since $\phi^{\prime}$ induces $|\phi^{\prime}|$ that sends
 points of $|\AM_0(X)_d|$ to closed points of
 $W_d^{\prime}$ that corresponds to
 $\Spec k\rightarrow W_d^{\prime}$, where $k$ is an algebraically
 closed field, any such morphism $\phi^{\prime}$ must send
 fiberwise isomorphic weak $\Delta$-collections on
 $(S\times{\Bbb P}^1)/S$ to the same element in
 $\Hom(S,W_d^{\prime})$.
This shows that $W_d$ corepresents $\AM_0(X)_d$.

\medskip

(b.1) and (b.2) together show that $\coprod_d W_d$ is the coarse
 moduli space for $\AM_0(X)$ and we conclude the proof.

\noindent\hspace{14cm}$\Box$

\bigskip

\bigskip

\section{The collapsing morphism.}
In this section, we re-run the proof of Jun Li of Lemma 2.7 in
 Mirror Principle II, with the A-twisted moduli stack $\AM_0(X)$ of
 Sec.\ 3 soldered into the discussion.
All the schemes in the discussion are over ${\Bbb C}$.

\bigskip

\begin{flushleft}
{\bf Background.}
\end{flushleft}
{\bf Fact 4.1 [rank $1$ sheaf].}
 ([Ha2]; also [Fr], [Huy-L], and [Od-S].) {\it
 Any rank $1$ torsion-free coherent sheaf on a locally factorial
  scheme $Y$ must be of the form ${\cal I}_Z\otimes {\cal L}$,
  where ${\cal I}_Z$ is the ideal sheaf of a subscheme $Z$ of codimension
  $\ge 2$ in $Y$ and ${\cal L}$ is a line bundle on $Y$.
 Such a decomposition is unique up to isomorphisms of
  ${\cal O}_Y$-modules.
} 

\bigskip

\noindent
{\bf Fact 4.2 [Hartogs extension theorem].} ([Ii].) {\it
 Let $Y$ be a Noetherian normal scheme and $Z$ be a closed subset
 of codimension $\ge 2$ in $Y$. Then
 $H^0(Y-Z,{\cal O}_Y)=H^0(Y,{\cal O}_Y)$.
} 

\bigskip

\noindent
{\bf Corollary 4.3 [determinant].} {\it
 Let ${\cal L}$ be a rank $1$ coherent sheaf on a locally
 factorial scheme $Y$. Then there exists a canonical morphism
 ${\cal L}\rightarrow\determinant{\cal L}$ of ${\cal O}_Y$-modules.
} 

\bigskip

\noindent
{\it Proof.}
 Recall the definition and the relations of $\,\determinant$ and
  $\,Div\,$ in [K-M].
 From the exact sequence
  $0\,\rightarrow\,\Tor{\cal L}\,\rightarrow\,{\cal L}\,
    \rightarrow\,{\cal L}/\!\mbox{\raisebox{-.4ex}{$\Tor{\cal L}$}}\,
    \rightarrow\,0\,$,
  one has
  $\determinant{\cal L}=\determinant(\Tor{\cal L})\otimes
    \determinant\left({\cal L}/
                \mbox{\raisebox{-.4ex}{$\Tor{\cal L}$}}\right)$.
 Since ${\cal L}/\mbox{\raisebox{-.4ex}{$\Tor{\cal L}$}}$ is
  torsion-free,
  ${\cal L}/\mbox{\raisebox{-.4ex}{$\Tor{\cal L}$}}
                   ={\cal I_Z}\otimes \widehat{\cal L}$ canonically,
  where
   ${\cal I}_Z$ is the ideal sheaf of the subscheme of codimension
     $\ge 2$ in $Y$ (from the flattening stratification, cf.\ [Mu3],
     of ${\cal L}/\mbox{\raisebox{-.4ex}{$\Tor{\cal L}$}}\,$)
     on which the fiber dimensions of
     ${\cal L}/\mbox{\raisebox{-.4ex}{$\Tor{\cal L}$}}$ jump up and
   $\widehat{L}$ is a line bundle on $Y$.
 These give rise to a sequence of canonical morphisms/identifications
  of ${\cal O}_Y$-modules$\,$:

  \vspace{-1ex}
  {\small
  $$
   {\cal L}\;
   \rightarrow\; {\cal L}/\mbox{\raisebox{-.4ex}{$\Tor{\cal L}$}}\;
   =\; {\cal I_Z}\otimes \widehat{\cal L}\;
   \rightarrow\; \widehat{\cal L}\;
   =\; \determinant\left({\cal I_Z}\otimes \widehat{\cal L}\right)\;
   \rightarrow\;
     \determinant\left({\cal I_Z}\otimes \widehat{\cal L}\right)
       \otimes \determinant\Tor{\cal L}\;
   \simeq \determinant{\cal L}\,,
  $$
  {\normalsize where}} 
  we have used the facts$\,$:
   (i)   $\determinant{\cal I}_Z={\cal O}_Y$,
   (ii)  $\determinant\Tor{\cal L}={\cal O}(\Div\Tor{\cal L})$
         and $\Div\Tor{\cal L}\ge 0$, and
   (iii) there are canonical inclusions
         ${\cal O}_Y\hookrightarrow {\cal O}_Y(D)
                                     \hookrightarrow {\cal K}_Y$
         for $D\ge 0$, where ${\cal K}_Y$ is the sheaf of total
         quotient rings of $Y$.
 The composition of this sequence of canonical morphisms gives
  the canonical morphism of ${\cal O}_Y$-modules
  ${\cal L}\rightarrow\determinant{\cal L}$ claimed.

\noindent\hspace{14cm}$\Box$

\bigskip

\noindent
{\bf Lemma 4.4 [push-pull of weak $\Delta$-collection].} {\it
 {\rm (1)}
  Let $f:Y\rightarrow Y^{\prime}$ be a dominant morphism that does
   not map an irreducible component of $Y$ to a point in $Y^{\prime}$,
   then the pull-back of a weak $\Delta$-collection on $Y^{\prime}$
   is a weak $\Delta$-collection on $Y$.

 \medskip

 \noindent {\rm (2)}
  Let $f:Y\rightarrow Y^{\prime}$ be a projective birational morphism
   between schemes of the same uniform dimension.
  Assume that $Y^{\prime}$ is irreducible and that $f$
   is an isomorphism outside a closed subscheme of codimension $\ge 2$
   in $Y^{\prime}$
   - in notation, $f|_{U}:U\stackrel{\sim}{\rightarrow}U^{\prime}$ -.
  Let $(L_{\rho}, u_{\rho}, c_m)_{\rho,m}$ be a weak
   $\Delta$-collection on $Y$ and
  define $L_{\rho}^{\prime}:= \determinant f_{\ast}L_{\rho}$.
  Then there exists a unique weak $\Delta$-collection
   $(L_{\rho}^{\prime},u_{\rho}^{\prime},c_m^{\prime})$
   on $Y^{\prime}$ that extends the weak $\Delta$-collection
   $(f|_U)_{\ast}(L_{\rho},u_{\rho},c_m)_{\rho,m}|_U$ on $U^{\prime}$.

 \medskip

 \noindent {\rm (3)}
  Let $f:Y=Y_0\cup Y_1\cup\cdots\,\rightarrow Y^{\prime}$
   be a projective morphism between schemes of the same uniform
   dimension that satisfies
   \begin{itemize}
    \item [{\rm (i)}]
     $Y^{\prime}$ is a Noetherian integral {\rm (}separated$\,${\rm )}
     scheme which is regular in codimension-$1$
      {\rm (}cf.\ {\rm [Ha1$\,$: II.6]}{\rm )},
    \item [{\rm (ii)}]
     the restriction $f:Y_0\rightarrow Y^{\prime}$ is an isomorphism
     outside a closed subscheme of codimension $\ge 2$ in $Y^{\prime}$
     - in notation, $f|_{U}:U\stackrel{\sim}{\rightarrow}U^{\prime}$ -,
     and
    \item [{\rm (iii)}]
     each $Y_i$, $i=1,\,\ldots$, is mapped to a codmension-$1$
     subscheme $D^{\prime}_i$ in $Y^{\prime}$, whose corresponding
     divisor is also denoted by $D^{\prime}_i$\,,
     {\rm (}i.e.\ $Y_i\rightarrow D^{\prime}_i$ is a flat family of
            curves{\rm )}.
   \end{itemize}
  Let $(L_{\rho}, u_{\rho}, c_m)_{\rho,m}$ be a weak
   $\Delta$-collection on $Y$ such that none of $u_{\rho}|_{Y_i}$
   are zero-sections, where $\rho\in\Delta(1)$ and $i=1,\,\ldots\,$,
   and let $L_{\rho}^{\prime}:= \determinant f_{\ast}L_{\rho}$.
  Then there exists a canonically constructed
   weak $\Delta$-collection
   $(L_{\rho}^{\prime}, u_{\rho}^{\prime}, c_m^{\prime})_{\rho,m}$
   on $Y^{\prime}$ that extends the weak $\Delta$-collection
   $(f|_{U-Y_1\cup\cdots})_{\ast}(L_{\rho},u_{\rho},c_m)_{\rho,m}
      |_{U-Y_1\cup\cdots}$
   on $U^{\prime}-D_1\cup\cdots\,$.
} 

\bigskip

\noindent
{\it Remark 4.5.}
 In Item (1), any condition that prevents mapping an irreducible
  component of $Y$ to the variety  $V(I)$ in ${\Bbb C}^{\Delta(1)}$
  will do. The condition stated here gives the kind of morphisms
  that appear in Jun Li's proof.
 Note that in Item (2) it is allowed that some $u_{\rho}$ are
  zero-sections while in Item (3) it is required that
  none of $u_{\rho}|_{Y_i}$, where $\rho\in\Delta(1)$ and
  $i=1,\,\ldots\,$, are zero-sections.

\bigskip

\noindent
{\it Proof.}
Statement (1) is clear and its counter statement for
 $\Delta$-collections over $\Spec{\Bbb C}$ is stated in [Cox2].

For Statement (2), observe that
$\otimes_{\rho}{L_{\rho}^{\prime}}^{\;\otimes\langle m,n_{\rho}\rangle}
 \simeq {\cal O}_{Y^{\prime}}$ as abstract
 ${\cal O}_{Y^{\prime}}$-modules since the former is an invertible
 ${\cal O}_{Y^{\prime}}$-module that is free outside a locus of
 codimension $\ge 2$ in $Y^{\prime}$.
The isomorphisms $(f|_U)_{\ast}c_m$ extend to unique isomorphisms
 $\otimes_{\rho}{L_{\rho}^{\prime}}^{\;\otimes\langle m,n_{\rho}\rangle}
  \simeq {\cal O}_{Y^{\prime}}$ by Hartogs extension theorem since,
 once fixing a trivialization of the two rank-$1$ globally free
 ${\cal O}_{Y^{\prime}}$-modules in question, $(f|_U)_{\ast}c_m$ is
 given by multiplication of a regular function.
The sections $u_{\rho}^{\prime}$ are given by the canonical morphism
 $H^0(Y,L_{\rho})\rightarrow H^0(Y^{\prime},L_{\rho}^{\prime})$
 arising from the combination of the definition of $f_{\ast}$ and
 the canonical morphism
 $f_{\ast}L_{\rho}
   \rightarrow L_{\rho}^{\prime}=\determinant f_{\ast}L_{\rho}\,$.
The cocycle conditions
 $c_{m_1}^{\prime}\otimes c_{m_2}^{\prime}=c_{m_1+m_2}^{\prime}$
 follow by continuity.

For Statement (3), let us first show that
 $\otimes_{\rho}{L_{\rho}^{\prime}}^{\;\otimes\langle m,n_{\rho}\rangle}
  \simeq {\cal O}_{Y^{\prime}}$
 as abstract ${\cal O}_{Y^{\prime}}$-modules.
By the assumption in Statement (3),
 $Y_i\rightarrow D_i$, $i=1,\,\ldots,\,$ are flat families of curves
 and the relative degree of $\reldeg_{D_i}(L_{\rho}|_{Y^i})$
 is well-defined. Moreover, recalling the definition of $\Div$,
 one concludes that
 $$
  L_{\rho}^{\prime}\; =\; \determinant f_{\ast}L_{\rho}\;
  =\;\determinant \left(\rule{0ex}{1.2em}
                     (f|_{Y_0})_{\ast}(L_{\rho}|_{Y_0})\right)\,
    \otimes\,
    {\cal O}_{Y^{\prime}}\left(
      \sum_{i=1,\,\ldots}\,\reldeg_{D_i}\,(L_{\rho}|_{Y_i})\,\cdot\,\,D_i
       \right)
 $$
 since all $u_{\rho}|_{Y_i}$ are non-zero sections.
By definition, the restriction
 $(L_{\rho}|_{Y_i}, u_{\rho}|_{Y_i}, c_m|_{Y_i})_{\rho, m}$
 of $(L_{\rho}, u_{\rho}, c_m)_{\rho,m}$ to each component $Y_i$ of $Y$
 is a weak $\Delta$-collection on $Y_i$. In particular,
 $c_m|_{Y_i}:
  \otimes_{\rho}\,L_{\rho}|_{Y_i}^{\;\langle m,n_{\rho}\rangle}
  \simeq {\cal O}_{Y_i}\,$ and
 $$
  \sum_{\rho\in \Delta(1)}\,\langle m,n_{\rho} \rangle\,\cdot\,
               \reldeg_{D_i}\,(L_{\rho}|_{Y_i}) \;=\; 0\,.
 $$
Furthermore, the restriction $f|_{Y_0}:Y_0\rightarrow Y^{\prime}$
 is in the situation of Statement (2) and one can define the weak
 $\Delta$-collection
 $(L_{\rho,0}^{\prime},
              u_{\rho,0}^{\prime}, c_{m,0}^{\prime})_{\rho,m}$
 on $Y^{\prime}$ as in Statement (2) as the push-forward of
 $(L_{\rho}, u_{\rho}, c_m)_{\rho,m}|_{Y_0}$ via $f|_{Y_0}$ .
It follows that
$$
 \otimes_{\rho\in\Delta(1)}\,
  {L_{\rho}^{\prime}}^{\,\otimes\langle m,n_{\rho}\rangle}\;
 =\; \otimes_{\rho\in\Delta(1)}\,{L_{\rho,0}^{\prime}}
              ^{\otimes\langle m,n_{\rho}\rangle}\;
  \stackrel{c_{m,0}^{\prime}}{\longrightarrow}\;
  {\cal O}_{Y^{\prime}}\,.
$$
This defines also the sought-for $c_m^{\prime}$.
 The sections $u_{\rho}^{\prime}$ and the cocycle conditions on
 $c_{m}^{\prime}$ follow by the same reasoning as in the case of
 Statement (2).
This concludes the proof.

\noindent\hspace{14cm}$\Box$

\bigskip

\bigskip

\begin{flushleft}
{\bf The collapsing morphism.}
\end{flushleft}

\noindent {\bf Proposition 4.6.} {\it
 Let 
  $X$ be a convex smooth toric variety and
  $M_d(X)$ be the moduli stack 
   $\overline{\cal M}_{0,0}({\Bbb P}^1\times X, (1,d))$ of genus $0$
   stable map into ${\Bbb P}^1\times X$ of degree $(1,d)$.
 Then there exists a natural morphism
  $\,\Upsilon:M_d(X)\rightarrow\AM_0(X)_d\,$ of stacks.
} 

\bigskip

\noindent
{\it Remark 4.7.} Composition of $\Upsilon$ with the morphism
 $\,\phi:\AM_0(X)_d\rightarrow W_d\,$ in the proof of
 Lemma 3.6 gives the morphism $\,\varphi:M_d(X)\rightarrow W_d\,$
 in [L-L-Y$\,$: II, Sec.\ 2.5, Lemma 2.7].

\bigskip

\noindent
{\it Proof of Proposition.}
We split the discussions to two cases.

\bigskip

\noindent
{\it Case $(1):$}
 \parbox[t]{13cm}{\it $M_d(X)$ is a compactification of the component
  of $\Hom({\Bbb P}^1,{\Bbb P}^1\times X)$ that corresponds to genus
  $0$ curves in ${\Bbb P}^1\times X$ of degree $(1,d)$.}

\medskip

\noindent
Let $(\Sch/\,{\Bbb C})$ be the category of Noetherian schemes of
 finite type over ${\Bbb C}$ and $\xi$ be an object in $M_d(X)(S)$
 given by
 $$
  \begin{array}{rcr}
   F:{\cal X}  & \longrightarrow  & S\times {\Bbb P}^1\times X \\
   & \searrow \hspace{1em} \swarrow & \\
   & S & .
  \end{array}
 $$
 Let $p_i$ (resp.\ $p_{ij}$) be the composition of $F$ with
 the projection of $S\times{\Bbb P}^1\times X$ to its $i$-th component
 (resp.\ the product of its $i$-th and $j$-th components).

Assume first that $(S,\xi)$ is an atlas of $M_d(X)$, then $S$ is smooth
 and outside a divisor $D_S$ of $S$ (i.e.\ the boundary locus of $S$)
 the defining family of stable maps over $S$ parameterizes morphisms
 from ${\Bbb P}^1$ into ${\Bbb P}^1\times X$.
The map
 $$
  p_{12}\;:\;{\cal X}\; \longrightarrow\;S\times {\Bbb P}^1
 $$
 is a projective birational morphism that is an isomorphism outside
 a locus of codimension $\ge 2$ in $S\times {\Bbb P}^1$.
Let $(L_\rho, z_\rho, c_m)_{\rho,m}$ be the universal
 $\Delta$-collection on $X$ and
 $$
  (L_{\rho,\xi}, u_{\rho,\xi}, c_{m,\xi})_{\rho,m}\;
  =\;p_3^{\ast}(L_\rho, z_\rho, c_m)_{\rho,m}\,.
 $$
 Then $(L_{\rho,\xi}, u_{\rho,\xi}, c_{m,\xi})_{\rho,m}$
  is a $\Delta$-collection on ${\cal X}$ over $S$.
The construction satisfies the base change property that,
 if $f:T\rightarrow S$ be a morphism of ${\Bbb C}$-schemes,
 then there is a canonical isomorphism of $\Delta$-collections
 $$
  (L_{\rho,f^{\ast}\xi}, u_{\rho,f^{\ast}\xi},
    c_{m,f^{\ast}\xi })_{\rho,m}\;
   \simeq\; (f\times \Id_{\,{\scriptsizeBbb P}^1})^{\ast}
             (L_{\rho,\xi}, u_{\rho,\xi}, c_{m,\xi})_{\rho,m}\,.
 $$
In particular, if one equips $M_d(X)$ with the \'{e}tale topology,
 then one obtains a $\Delta$-collection on the stack $M_d(X)$
 by considering the \'{e}tale morphisms among atlases.

Let ${\cal L}_{\rho,\xi}={p_{12}}_{\ast}L_{\rho,\xi}$.
Since $p_{12}$ is an isomorphism over the complement $U^{\prime}$
 of a codimension $\ge 2$ locus in $S\times{\Bbb P}^1$,
 by Lemma 4.4\,(2) there exists a unique weak $\Delta$-collection
 of the form
 $(\determinant{\cal L}_{\rho,\xi},
   \sigma_{\rho,\xi}, c_{m,\xi}^{\;\prime})$
 on $S\times{\Bbb P}^1$ (over $\Spec{\Bbb C}$) such that the
 restriction of ${p_{12}}_{\ast}$ over $U^{\prime}$ is an isomorphism
 of weak $\Delta$-collections on $U^{\prime}$ (over $\Spec{\Bbb C}$).
Since
 each fiber ${\Bbb P}^1$ of $S\times {\Bbb P}^1$ over $S$ comes from
  pinching rational subcurves of the corresponding fiber of ${\cal X}$
  over $S$ and
 the restriction of
  $(\determinant{\cal L}_{\rho,\xi},
               \sigma_{\rho,\xi}, c_{m,\xi}^{\;\prime})$
  to a fiber ${\Bbb P}^1$ defines a morphism from ${\Bbb P}^1$ to
  $X$ (of possibly lower multi-degrees),
$(\determinant{\cal L}_{\rho,\xi},
                \sigma_{\rho,\xi}, c_{m,\xi}^{\;\prime})$
 must satisfy the nonvanishing condition of
 Definition 2.1.3   
 when restricted to each fiber ${\Bbb P}^1$ of $S\times{\Bbb P}^1$
 over $S$.
Consequently,
 $(\determinant{\cal L}_{\rho,\xi},
               \sigma_{\rho,\xi}, c_{m,\xi}^{\;\prime})$
 is a weak $\Delta$-collection on $S\times{\Bbb P}^1$ over $S$ as well
and one obtains a map
 $$
  \Omega\,:\,\{\,\mbox{atlases $(S,\xi)$ of $M_d(X)$}\,\}\;
                     \longrightarrow\; \AM_0(X)_d
 $$
 that commutes with the \'{e}tale base change among atlases of
 $M_d(X)$.

Fix now an atlas $(T,\xi_T)$ for $M_d(X)$ and let $\xi\in M_d(X)(S)$
 for a general $S\in (\Sch/\,{\Bbb C})$.
Since $M_d(X)$ is a smooth Deligne-Mumford stack, the pair
 $(\xi_T, \xi)$ determines a commutative diagram
 $$
  \begin{array}{ccc}
   S^{\prime}:=\Isombf(\xi_T,\xi)
     & \stackrel{\alpha}{\longrightarrow}
     & S  \\
   \hspace{2ex}\downarrow\mbox{\scriptsize $\beta$} & & \downarrow \\
   T & \longrightarrow & M_d(X)\,,
  \end{array}
 $$
 where $\alpha$ is \'{e}tale and surjective.
The canonical isomorphism $\alpha^{\ast}\xi\simeq \beta^{\ast}\xi_T$
 induces a canonical isomorphism
 $\alpha^{\ast}(L_{\rho,\xi},u_{\rho,\xi},c_{m,\xi})_{\rho,m}
  \simeq \beta^{\ast}(L_{\rho,\xi_T},u_{\rho,\xi_T},
                                            c_{m,\xi_T})_{\rho,m}$.
The weak $\Delta$-collection $\beta^{\ast}\Omega(T,\xi_T)$ on
 $(S^{\prime}\times{\Bbb P}^1)/S^{\prime}$ is a descent datum with
 respect to $\alpha$ and hence descends to a weak $\Delta$-collection
 on $(S\times{\Bbb P}^1)/S$.
One can check that different choices of atlases $(T,\xi_T)$ give
 rise to the same descent on $(S\times{\Bbb P}^1)/S$, thus one obtains
 a well-defined morphism of stacks from $M_d(X)$ to
 $\AM_0(X)_d$.
This concludes the proof for Case (1).

\bigskip

\noindent
{\it Case $(2):$}
 \parbox[t]{13cm}{\it General $M_d(X)$. }

\medskip

\noindent
Again, let $S$ be an atlas of $M_d(X)$, which is smooth.
 Then there is a stratification of $S$ labelled by the dual
 graphs of the prestable domain curves of stable maps in question.
We shall assume that the graph for the maximal stratum is not a point,
 i.e.\ we are not in Case (1).
Then the projective morphism $p_{12}\,$ are now in the situation of
 Statement (3) of Lemma 4.4.
Convexity of $X$ implies that the restriction of $u_{\rho,\xi}$,
 as defined analogous to the discussion in Case (1), to each component
 of ${\cal X}$ is not a zero-section.
 The proposition now follows from Lemma 4.4\,(3) and the same argument
 as in Case (1) above.

\noindent\hspace{14cm}$\Box$

\bigskip

We conclude the notes with three themes along the line for
 further study.

\bigskip

\noindent
{\bf Theme 1.}
 Further properties and details of the A-twisted moduli stack
 $\AM_g(X)$.

\bigskip

\noindent
{\bf Theme 2.}
 Construction of natural morphisms between the moduli stack of
  stable maps and the A-twisted moduli stack that generalize
  the construction in [L-L-Y].

\bigskip

\noindent
{\bf Theme 3.}
 Generalization of the twisted moduli stack $\AM_g(X)$ to
 the case of open strings, e.g.\ [G-J-S].

\bigskip

\bigskip

\begin{flushleft}
{\large\bf  
  Appendix. Witten's gauged linear sigma models for mathematicians.}
\end{flushleft}
Witten's gauged linear sigma model (GLSM) [Wi1] is one of the universal
 frameworks or structures that lie behind stringy dualities (e.g.\ [Gre]).
A mathematical review of the related part of [Wi1] (cf.\ also [M-P])
 to the current work is given in this subsection.

\bigskip

\noindent $\bullet$
{\bf Introduction to the superland.}
[Po2$\,$: vol.\ II.\ Appendix B] (resp.\ [Fr]) gives a concise
 introduction of {\it spinor representations}, {\it supersymmetry} (SUSY),
 {\it supermultiplets},
 and {\it superfields} and their {\it component fields}
 from a string theorist's (resp.\ mathematician's) aspect.
A formulation of {\it superspaces/manifolds/schemes} that is close
 in spirit to Grothendieck's formulation of algebraic geometry
 is given in [Ma$\,$: Chapter 4].
This formulation provides a geometry behind the standard text [W-B]
 on supersymmetry.
{\it K\"{a}hler differentials} and {\it tangent vectors} can be
 defined as in [Ha1].
{\it Fermionic integration} is discussed in [We$\,$: Sec.\ 26.6]
 and [W-B$\,$: IX], whose mathematical formulation {\it Berezin integral}
 is discussed in [Fr] and [Ma].
{\it R-symmetry} is discussed in [Fr$\,$: Lecture 3] and [We].
{\it Central extensions} of a supersymmetry algebra and
 its {\it BPS representations} are discussed in [Fr], [Po2$\,$: vol II],
 and [We].
{\it Super linear algebra},
 in particular the {\it parity change functor} $\prod$,
 is discussed in [Fr$\,$: Lectures 1 and 2] and [Ma$\,$: Chaper 3].
See also [Arg] and [DEFJKMMW].

\bigskip

\noindent $\bullet$
{\bf Supermanifolds and line bundles.}
For the purpose of this article we have reduced the role played
 by the parity change functor $\prod$ in the description as much
 as possible.
\begin{itemize}
 \item [(1)] ([Ma$\,$: Sec.\ 4.1]; also [Ha1].)
  A {\it supermanifold} $X=(X, {\cal O}_X)$ (in smooth, analytic,
   or algebraic category) is a ${\Bbb Z}/2{\Bbb Z}$-graded ringed
   topological space $(X,{\cal O}_X)$ such that
   \begin{itemize}
    \item [$(a)$]
     The stalk ${\cal O}_{X,x}$ of ${\cal O}_X$ at any point $x\in X$
     is a local ring.

    \item [$(b)$]
     $X$ is covered by a collection of open sets
      $\{U_{\alpha}\}_{\alpha\in I}$ such that
      each $(U_{\alpha},{\cal O}_M|_{U_{\alpha}})$ is isomorphic to
      $(\,U^{\,0}_{\alpha}\,,\,
         \Sym^{\bullet}_{\,{\cal O}_{U_{\alpha}^0}}(
                                         \prod{\cal E}_{\alpha})\,)$,
      where
       $U^{\,0}_{\alpha}$ is an ordinary manifold
        (in the corresponding category),
       ${\cal E}_{\alpha}$ is an ordinary locally free coherent
        ${\cal O}_{U^{\,0}_{\alpha}}$-module, and
       $\prod{\cal E}_{\alpha}$
        is the ${\cal O}_{U^{\,0}_{\alpha}}$-module ${\cal E}_{\alpha}$
        with odd parity.

    \item [$(c)$]
     $X_{\rm rd}$ is a manifold (in the corresponding category),
     cf.\ Remark below.
   \end{itemize}
   {\it Remark.}
    Let
     ${\cal O}_X={\cal O}_X^{(0)}\oplus{\cal O}_X^{(1)}$
      be the decomposition of ${\cal O}_X$ into the even
      (i.e.\ grade $0$) and the odd (i.e.\ grade $1$) component and
     $J_X := {\cal O}_X^{(1)}+({\cal O}_X^{(1)})^2$
      (the ideal of ``superfuzz", cf.\ [Fr$\,$: Lecture 1]).
    Then $X_{\rm rd}$ is by definition the submanifold of $X$
     associated to $J_X$.
    Note also that
     $\Sym^{\bullet}(\prod{\cal E}_{\alpha})
          \simeq \bigwedge^{\bullet}{\cal E}_{\alpha}$
     as ${\cal O}_{U^{\,0}_{\alpha}}$-modules with all the parities
     after tensor products erased.
    $X$ is called {\it decomposable} if in Condition (b) one can
     choose $U_{\alpha}=X$ for some $\alpha$.
    In this case there is a surjective affine morphism
     $X\rightarrow X_{\rm rd}$ such that the composition
     $X_{\rm rd}\hookrightarrow X\rightarrow X_{\rm rd}$
     is the identity map.

 \item [(2)]
  A {\it line bundle} ${\cal L}$ on $X$ is a locally free rank $1$
   ${\cal O}_X$-module.
  Associated to ${\cal L}$ is a finite filtration of
   ${\cal O}_X$-modules$\,$:
   ${\cal L}\,\supset\, {\cal L}\cdot J_X\,
     \supset\,{\cal L}\cdot J_X^{\,2}\,\supset\,\cdots\,\supset\,0\,$.
  Global sections of the associated graded object
   $\Gr{\cal L}
    :=\oplus_{\,i}\,({\cal L}\cdot J_X^{\,i}/{\cal L}\cdot J_X^{\,i+1})$
   are called {\it component sections} of ${\cal L}$.
  The restriction ${\cal L}_{\rm rd}$ of ${\cal L}$ to $X_{\rm rd}$
   is a usual line bundle on $X_{\rm rd}$.
  When $X$ is decomposable,
   $\Gr{\cal L}\simeq{\cal O}_M\otimes{\cal L}_{\rm rd}$
   as ${\cal O}_{X_{\rm rd}}$-modules with all the parities erased.
  One may define also the Picard group $\Pic(X)$ of $X$.
\end{itemize}

\bigskip

\noindent $\bullet$
{\bf $N\,$: the count of minimal collections.}
([Fr$\,$: Lecture 3] and [Po2$\,$: vol.\ II, Appendix B].)
 The real dimension of a minimal real spin representation at
  $d$-dimensional Minkowski space is given by
  $$
   \mbox{
    \begin{tabular}{c|cccccccccccc}
     $d$  & 1 & 2 & 3 & 4 & 5 & 6 & 7 & 8 & 9 & 10 & 11 & 12 \\ \hline
     $\dimm_{\scriptsizeBbb R}$
          & 1 & 1 & 2 & 4 & 8 & 8 & 16 & 16 & 16 & 16 & 32 & 64
    \end{tabular}
    }\,.
  $$
 In even dimensions, there are two such irreducible representations,
  distinguished by {\it left} and ${\it right}$.
 The $N$ that appears in every SUSY literatures counts the number of
  collections of the odd generators of a SUSY algebra with each
  collection in a minimal spinor representation of the Lorentz
  subalgebra of the SUSY algebra.

\bigskip

\noindent
{\bf Example A.1 [$d=4, N=2$].}
 The complexified $d=4,\, N=2$ SUSY algebra
  (as in the Seiberg-Witten theory) contains $8$ odd generators
  in collections of $4$.
 Each collection of odd generators spans an irreducible represenation,
  ${\bf 2}$ or ${\bf 2^{\prime}}$, of $\Spin(1,3)\simeq \SL(2,{\Bbb C})$.
 For $d=4$, though there are two different minimal spinor representations,
   they give the same complexification ${\bf 2}+{\bf 2^{\prime}}$.
 In physics literature SUSY algebras are ususally complexified; thus
  it is not necessary to distinguish whether it is ${\bf 2}$ or
  ${\bf 2^{\prime}}$ that appears in the SUSY algebra at $d=4$.
 In contrast, at $d=2$, complexifications of ${\bf 1}$ and
  ${\bf 1}^{\prime}$ give inequivalent representations of
  $\Spin(1,1)$ and the distinction of left and right is necessary.
 E.g.\ $N=(1,1)$ and $N=(0,2)$ label different complexified SUSY
  algebras at $d=2$ with $2$ odd generators. The distinction is also
  needed at $d=10$, cf.\ the mod-$8$ periodicity of many properties
  of spinor representations.

\bigskip

\noindent $\bullet$
{\bf Physical supermanifolds.}
For simplicity and sufficiency of this paper, we shall assume that
 the supermanifold $X$ is decomposable,
 i.e.\ ${\cal O}_X=\Sym^{\bullet}_{{\cal O}_{X_{\rm rd}}}(\prod{\cal E})$
 for some locally free ${\cal O}_{X_{\rm rd}}$-module ${\cal E}$.
To link $X$ with supersymmetry from physics, it is then required
 that $X_{\rm rd}$ is a Lorentzian manifold and ${\cal E}$ is
 a spinor bundle on $X_{\rm rd}$.
{\it Superfields} on $X$ are defined to be global sections of locally
 free sheaves, e.g.\
 ${\cal O}_X=\Sym^{\bullet}_{{\cal O}_{X_{\rm rd}}}(\prod{\cal E})$,
 on $X$.

\bigskip

\noindent
{\bf Example A.2 [$d=4$, $N=1$].}
(1) {\it The supergeometry.}
$X_{\rm rd}=$ the Minkowski space (with coordinates $x=(x^0,x^1,x^2,x^3)$)
  equipped with the standard metric of signature $(-1,1,1,1)$ and
 ${\cal E}=({\cal O}_{X_{\rm rd}}\otimes{\bf 2})_{\,\scriptsizeBbb C}
   =({\cal O}_{X_{\rm rd}})_{\,\scriptsizeBbb C}
         \otimes_{\scriptsizeBbb C}({\bf 2}+{\bf 2^{\prime}})$,
 the complexified spinor bundle on $X_{\rm rd}$.
Fix a set of (anticommuting) generators $\theta^{\alpha}$,
 $\overline{\theta}^{\dot{\alpha}}$, $\alpha=1,2$, for $\prod{\cal E}$
 as an $({\cal O}_{X_{\rm rd}})_{\,\scriptsizeBbb C}$-module.
Recall the decomposition
 $({\bf 2}+{\bf 2^{\prime}})\bigwedge({\bf 2}+{\bf 2^{\prime}})
   = {\bf 1} + {\bf 1} + {\bf 4}$,
 where
  ${\bf 1}$ is the (complexified) $1$-dimensional trivial
   representation and
  ${\bf 4}$ is the (complexified) vector represenattion of $\SO(1,3)$,
 the Pauli matrices $\sigma^m$, $m=0,1,2,3$ and the $\varepsilon$
 matrices (cf.\ [W-B$\,$: Appendix B] and [Fr$\,$: Lecture 3]).
Then a superfield from $({\cal O}_X)_{\scriptsizeBbb C}$ can be
 expressed as (cf.\ [W-B$\,$: Appendix A] for summation conventions)
 \begin{eqnarray*}
  \lefteqn{
   F(x,\theta,\overline{\theta})\;
    =\; f(x) + \theta\phi(x) + \overline{\theta}\overline{\chi}(x)
       + \theta\theta m(x) + \overline{\theta}\overline{\theta} n(x)
       + \theta\sigma^m\overline{\theta}v_m(x)    }\\[.6ex]
   & & \hspace{18em}
       + \theta\theta\overline{\theta}\overline{\lambda}(x)
       + \overline{\theta}\overline{\theta}\theta\psi(x)
       + \theta\theta\overline{\theta}\overline{\theta} d(x)
 \end{eqnarray*}
 with the component fields from representations of $\SO(1,3)\,$:
 $$
  \begin{array}{ccccccccc}
    1                    & \theta            & \overline{\theta}
        & \theta\theta    & \overline{\theta}\overline{\theta}
        & \theta\sigma^m\overline{\theta}
        & \theta\theta\overline{\theta}
        & \overline{\theta}\overline{\theta}\theta
        & \theta\theta\overline{\theta}\overline{\theta} \\[.6ex]
     f(x)  & \phi(x)  & \overline{\chi}(x)  & m(x)  & n(x)  & v_m(x)
           & \overline{\lambda}(x) & \psi(x)  & d(x) \\[.6ex]
     {\bf 1}        & {\bf 2}     & {\bf 2^{\prime}}
        & {\bf 1}   & {\bf 1}
        & {\bf 4}   & {\bf 2^{\prime}}     & {\bf 2}
        & {\bf 1}   \\
   \end{array}
  $$
  (cf.\ [W-B$\,$: Eq.(4.9)]).
The $d=4$, $N=1$ SUSY algebra can be realized as an algebra of
 (differential) operators acting an $F(x,\theta,\overline{\theta})$.
In particular the four odd generators are realized as$\,$:
 (This is what physicists call {\it SUSY generators} of the SUSY algebra.)
 $$
  Q_{\alpha}\;
   =\; \frac{\partial}{\partial\theta^{\alpha}}
        - \sqrt{-1}\,\sigma^m_{\alpha\dot{\alpha}}\,
                       \overline{\theta}^{\dot{\alpha}}
          \frac{\partial}{\partial x_m}
   \hspace{1em}\mbox{and}\hspace{1em}
  \overline{Q}_{\dot{\alpha}}\;
   =\; -\,\frac{\partial}{\partial\overline{\theta}^{\dot{\alpha}}}
        + \sqrt{-1}\,\sigma^m_{\alpha\dot{\alpha}}\,\theta^{\alpha}
          \frac{\partial}{\partial x_m}\,,
   \hspace{1ex}\alpha,\dot{\alpha}=1,\,2\,,
 $$
 ([Wi1$\,$: Eq.(2.1)] and [W-B$\,$: Eq.(4.4)]).

\medskip

\noindent
(2) {\it Chiral superfields and chiral multiplets.}
A superfield $\Phi$ (resp.\ $\overline{\Phi}$) that satisfies
$$
 \overline{D}_{\dot{\alpha}}\Phi\;=\;0 \hspace{1em}
  (\,\mbox{resp.}\hspace{1ex} D_{\alpha}\overline{\Phi}\;=\;0\,)\,,
$$
where
$$
 D_{\alpha}\;
  =\; \frac{\partial}{\partial\theta^{\alpha}}
       + \sqrt{-1}\,\sigma^m_{\alpha\dot{\alpha}}\,
                      \overline{\theta}^{\dot{\alpha}}
         \frac{\partial}{\partial x_m}
  \hspace{1em}\mbox{and}\hspace{1em}
 \overline{D}_{\dot{\alpha}}\;
  =\; -\,\frac{\partial}{\partial\overline{\theta}^{\dot{\alpha}}}
       - \sqrt{-1}\,\sigma^m_{\alpha\dot{\alpha}}\,\theta^{\alpha}\,
         \frac{\partial}{\partial x_m}\,,
$$
is called {\it chiral superfield} (resp.\ {\it antichiral superfield}).
In terms of $y^m=x^m+\sqrt{-1}\theta\sigma^m\overline{\theta}$
 in the coordinate ring of $X$, such $\Phi$ (resp.\ $\overline{\Phi}$)
 can be expressed as
 $$
  \Phi(y,\theta)\;=\;\phi(y)+\sqrt{2}\,
                          \theta\psi(y)+\theta\theta F(y)
 $$
(resp.\
 $$
  \overline{\Phi}(\overline{y},\overline{\theta})\;
   =\; \overline{\phi}(\overline{y})
      + \sqrt{2}\,\overline{\theta}\overline{\psi}(\overline{y})
      + \overline{\theta}\overline{\theta}\overline{F}(\overline{y})\,,
 $$
) in component fields.
The part $(\phi,\psi)$ is a section of the vector bundle associated to
 a $d=4$, $N=1$ {\it chiral multiplet} representation
 (cf. [Fr$\,$: Lecture 5, Table 7]).
When the Lagrangian for $d=4, N=1$ SUSY quantum field theory (SQFT)
 is considered, the equation of motion for $(\phi,\psi)$ will involve
 differential operators while that for $F$ will be purely algebraic.
 We say that $\phi$ and $\psi$ are {\it dynamical} component fields and
 $F$ {\it auxiliary} component field in the chiral multiplet $\Phi$.

\medskip

\noindent
(3) {\it Vector superfields and vector multiplets.}
 A superfield $V$ that satisfies the reality condition
  $$
   V\;=\; V^{\dagger}\,,
  $$
  where $V^{\dagger}$ is the Hermitian conjugate of $V$
  ([W-B$\,$: Appendix A]), is called a {\it vector superfield}.
 In the {\it Wess-Zumino gauge}, its component field expansion is
  $$
   V\; =\; -\theta\sigma^m\overline{\theta}v_m
    +\sqrt{-1}\,\theta\theta\overline{\theta}\overline{\lambda}
    -\sqrt{-1}\,\overline{\theta}\overline{\theta}\theta\lambda
    +\frac{1}{2}\,\theta\theta\overline{\theta}\overline{\theta}\,D\,,
  $$
 ([W-B$\,$: Eq.(6.6)] and [Wi1$\,$: Eq.(2.11)]).
 The dynamical components $(\lambda, \overline{\lambda}, v_m)$
  is a section of the vector bundle associated to the $d=4$, $N=1$
  {\it massless vector multiplet} representation
  (cf. [Fr$\,$: Lecture 5, Table 7])
  while $D$ is an auxiliary component, which plays an important role
  in defining the vacuum manifolds in each phase of a gauged linear
  sigma model [Wi1].

\bigskip

\noindent $\bullet$
{\bf Dimensional reduction $(d=4, N=1)\Rightarrow(d=2, N=(2,2))$
     and R-symmetry.} (Cf.\ [DEFJKMMW], [H-V], [We], and [Wi1].)
\begin{itemize}
 \item [(1)]
  The {\it $d=4,\, N=1$ SUSY algebra} is given by generators with
   (anti-)commutation relations$\,$: ($(\eta_{mn})=\Diag(-1,1,1,1)$.)

  \vspace{-1.6em}
  \item []
  {\footnotesize
  $$
   \begin{array}{lr}
    [L_{mn},L_{m^{\prime}n^{\prime}}]\;
     =\;   \eta_{nm^{\prime}}L_{mn^{\prime}}\,
        -\,\eta_{mn^{\prime}}L_{nn^{\prime}}\,
        -\,\eta_{n^{\prime}m}L_{m^{\prime}n}\,
        +\,\eta_{n^{\prime}n}L_{m^{\prime}m}\,,
           \hspace{-10em}& (\,\mbox{Lorentz algebra}\,)      \\[.6ex]
    [L_{mn},P_{m^{\prime}}]\;
      =\;\eta_{m^{\prime}n}P_m - \eta_{m^{\prime}m}P_n
                 & (\,\mbox{vector representation}\,)   \\[.6ex]
    [L_{mn},Q_{\alpha}]\;
      =\; (\sigma_{mn})_{\alpha}^{\;\;\beta}\,Q_{\beta}\,,
       \hspace{1em}
    [L_{mn},Q_{\dot{\alpha}}]\;
      =\; Q_{\dot{\beta}}\,
           (\overline{\sigma}_{mn})_{\;\;\dot{\alpha}}^{\dot{\beta}}\,,
                 & (\,\mbox{spinor representation}\,)   \\[.6ex]
    \{Q_{\alpha},\overline{Q}_{\dot{\alpha}}\}\;
      =\;2\sigma_{\alpha\dot{\alpha}}^m\,P_m
                 & (\,\mbox{Clifford-type algebra}\,)   \\[.6ex]
    \{Q_{\alpha},Q_{\beta}\}\;
      =\;\{\overline{Q}_{\dot{\alpha}},\overline{Q}_{\dot{\beta}}\}\;
      =\;0\,                                            \\[.6ex]
    [J,Q_{\alpha}]\; =\; Q_{\alpha}\,,\hspace{1em}
     [J,\overline{Q}_{\dot{\alpha}}]\;=\; -\,\overline{Q}_{\dot{\alpha}}
                    & (\,\mbox{R-symmetry $U(1)$}\,)    \\[.6ex]
    [P_m,P_n]\;=\; [J,P_m]\; =\; [J,L_{mn}]=0\,.        \\[.6ex]
    (\,m,\, n,\, m^{\prime},\, n^{\prime} = 0, 1, 2, 3\,;\;
     \alpha,\, \beta\, =1,2\,;\;
     \dot{\alpha},\, \dot{\beta} = \dot{1}, \,\dot{2}\,,
      \hspace{1ex}\mbox{[W-B$\,$: Eq.(A.14)]}\,.\,) \hspace{-3em}&
   \end{array}
  $$
  {\normalsize It}} 
  is customary to call $Q_{\alpha}$, $\overline{Q}_{\dot{\alpha}}$
  the {\it SUSY generators} of the SUSY algebra.

 \item [(2)]
  The {\it dimensional reduction} of the $d=4,\, N=1$ SUSY algebra
   to $d=2$ is obtained by considering the subalgebra that leaves
   a specified ${\Bbb R}^{1+1}$ subspace, e.g.\
   the $(x^0,x^3)$-coordinate plane, in ${\Bbb R}^{1+3}$ invariant.
  This corresponds to setting the extra conditions
   $$
    P_1\; =\; P_2\;
    =\; L_{01}\; =\; L_{02}\; =\; L_{13}\; =\; L_{23}\; =\; 0
   $$
   to the $d=4,\, N=1$ SUSY algebra since these generators generate
   Lorentz transformations that do not leave the $(x^0,x^3)$-coordinate
   plane invariant.
  The resulting algebra is the {\it $d=2,\, N=(2,2)$ SUSY algebra}.
  Its generators with renamings are

  \vspace{-1.2em}
  \item[]
   {\footnotesize
   $$
    \begin{array}{l}
     L:=L_{03},\hspace{2ex}  H:=-P_0,\hspace{2ex}  P:= P_3,\hspace{2ex}
     Q_- := Q_1,\hspace{2ex}     Q_+ := Q_2,\hspace{2ex}
      \overline{Q}_- := \overline{Q}_{\dot{1}},\hspace{2ex}
      \overline{Q}_+ := \overline{Q}_{\dot{2}},     \\[.6ex]
     J_1 := J,\hspace{2ex}
     J_2 := -\,2\sqrt{-1}\,L_{12}\,,
    \end{array}
   $$
   {\normalsize with}} 
   commutation relations$\,$:
   {\footnotesize
   $$
    \begin{array}{l}
     [L,H]\;=\;  -\,P\,, \hspace{2em}  [L,P]\;=\; -\,H\,,   \\[.6ex]
     [L,Q_+]\;=\;     \frac{1}{2}Q_+\,, \hspace{1.6em}
      [L,Q_-]\;=\; -\,\frac{1}{2}Q_-\,,   \hspace{1.6em}
      [L,\overline{Q}_+]\; =\; \frac{1}{2}\overline{Q}_+\,,\hspace{1.6em}
      [L,\overline{Q}_-]\; =\; -\,\frac{1}{2}\overline{Q}_-\,, \\[.6ex]
     \{Q_+,\overline{Q}_+\}\;=\; 2\,(H-P)\,, \hspace{2em}
      \{Q_-,\overline{Q}_-\}\;=\; 2\,(H+P)\,,    \\[.6ex]
     [J_1,Q_+]\;=\;   Q_+\,,       \hspace{2em}
      [J_1,Q_-]\;=\;   Q_-\,,       \hspace{2em}
      [J_1,\overline{Q}_+]\;=\; -\,\overline{Q}_+\,, \hspace{2em}
      [J_1,\overline{Q}_-]\;=\; -\,\overline{Q}_-           \\[.6ex]
     [J_2,Q_+]\;=\;    Q_+\,,       \hspace{2em}
      [J_2,Q_-]\;=\; -\,Q_-\,,       \hspace{2em}
      [J_2,\overline{Q}_+]\;=\; -\,\overline{Q}_+\,, \hspace{2em}
      [J_2,\overline{Q}_-]\;=\;    \overline{Q}_-           \\[.6ex]
     Q_+^2\; =\;Q_-^2\; =\;\overline{Q}_+^2\; =\;\overline{Q}_-^2\;
      =\;\{Q_-,Q_+\}\; =\;\{\overline{Q}_-,\overline{Q}_+\}\;
      =\; \{Q_+,\overline{Q}_-\}\; =\;\{Q_-,\overline{Q}_+\}\;
      =\;0\,,  \\[.6ex]
     [H,P]\;
      =\;[J_1,H]\;=\;[J_1,P]\;=\;[J_1,L]\;
      =\;[J_2,H]\;=\;[J_2,P]\;=\;[J_2,L]\;
      =\;[J_1, J_2]\;=\;0\,.
    \end{array}
   $$
   }
   $Q_-$ and $\overline{Q}_-$ (resp.\ $Q_+$ and $\overline{Q}_+$)
    are the $d=2,\, N=(2,0)$ (resp.\ $N=(0,2)$) SUSY generators and
   the Lorentz generator $L_{12}$ in the original algebra has now
    become the second R-symmetry generator $J_2$ of the new SUSY algebra.

  \vspace{-.8ex}
  \item []
  \hspace{1.6em}Define
   $$
    \begin{array}{ll}
     J_L\; =\; \frac{1}{2}\,(J_2-J_1)\,,
      & (\,\mbox{left-moving R-symmetry generator}\,) \\[.6ex]
     J_R\; =\; \frac{1}{2}\,(J_2+J_1)\,.
      & (\,\mbox{right-moving R-symmetry generator}\,)
    \end{array}
   $$
   Then
   {\footnotesize
   $$
    \begin{array}{llll}
     [J_L,Q_-]\;=\; -\,Q_-\,,
       & [J_L,\overline{Q}_-]\;=\; \overline{Q}_-\,,
       & [J_L,Q_+]\;=\; 0\,,
       & [J_L,\overline{Q}_+]\;=\; 0\,        \\[.6ex]
     [J_R,Q_-]\;=\; 0\,,
       & [J_R,\overline{Q}_-]\;=\; 0\,.
       & [J_R,Q_+]\;=\; Q_+\,,
       & [J_R,\overline{Q}_+]\;=\; -\,\overline{Q}_+\,,
    \end{array}
   $$
   } 

 \item [(3)]
  The dimensional reduction of a superfield on $d=4$ Minkowski
   space-time to $d=2$ superfields is obtained by setting two spatial
   directions, say $x^1$ and $x^2$, to be constant and take fields
   to depend only on $x^0$ and $x^3$.
  Recall Example A.2. 
  Then the rule of conversion of fields from $d=4,\,N=1$ to
   $d=2,\,N=(2,2)\,$ are given by$\,$:

 \vspace{-3ex}
 \item []
 {\footnotesize
   $$
    \begin{array}{lr}
     y^0 :=x^0\,,\hspace{3ex}  y^1 :=x^3\,,
                  & \mbox{(for space-time coordinates)}   \\[.6ex]
     \sigma:= (v_1-\sqrt{-1}v_2)/\sqrt{2}\,,\hspace{3ex}
     \overline{\sigma}:=(v_1+\sqrt{-1}v_2)/\sqrt{2}
                  &
                      \mbox{\parbox[t]{38ex}{(reduced vector components
                        in $d=4$ \newline
                        $\Rightarrow$ complex scalors in $d=2$)}}\\[.6ex]
     \left\{
      \begin{array}{cc}
       (\psi^-,\psi^+) :=(\psi^1,\psi^2)\,,
         & (\psi_-,\psi_+) :=(\psi_1,\psi_2)\,, \\
       (\overline{\psi}^-,\overline{\psi}^+)
          :=(\psi^{\dot{1}},\psi^{\dot{2}})\,,
         & (\overline{\psi}^-,\overline{\psi}^+)
          :=(\psi^{\dot{1}},\psi^{\dot{2}})\,.
      \end{array}
     \right.
      & \mbox{(for spinorial components)}
    \end{array}
   $$
  {\normalsize The}} 
  $d=2\,N=(2,2)$ chiral superfields and vector superfields
  can be obtained from $d=4,\, N=1$ ones via these conversions,
  cf. Example A.2 
  and the next item.
\end{itemize}

\bigskip

\noindent $\bullet$
{\bf Gauged linear sigma models.} ([Wi1$\,$: Sec.\ 2].)
Given $d=2,\,N=(2,2)$ chiral superfields
($\Phi_i$ here is the $\Phi_0$ in [Wi1$\,$: Eq.(2.13)])
 $$
  \Phi_i\;=\; \phi_i + \sqrt{2}\theta^+\psi_{i,+}
              + \sqrt{2}\theta^-\psi_{i,-}
              + \theta^2 F_i\,,
   \hspace{1em}i=1,\,\ldots,\,n\,,
 $$
 in $d=2,\,N=(2,2)$ chiral
 coordinates ($\Phi_i$ here is the $\Phi_0$ in [Wi1$\,$: Eq.(2.13)])
 $$
  (\,y^0-\sqrt{-1}
       (\theta^-\overline{\theta}^- +\theta^+\overline{\theta}^+)\,,\,
      y^1+\sqrt{-1}
       (\theta^-\overline{\theta}^- -\theta^+\overline{\theta}^+)\,)\,,
 $$
 vector superfields, in Wess-Zumino gauge, (cf.\ connections)

 \vspace{-1ex}
 {\footnotesize
  \begin{eqnarray*}
   \lefteqn{
    V_a\;
     =\; -\sqrt{2}\,
           (\theta^-\overline{\theta}^+\sigma_a
             + \theta^+\overline{\theta}^-\overline{\sigma}_a)
      + (\theta^-\overline{\theta}^-
             + \theta^+\overline{\theta}^+) v_{a,0}
      - (\theta^-\overline{\theta}^-
             - \theta^+\overline{\theta}^+) v_{a,1}}\\[.6ex]
    & & \hspace{2em}
     + \sqrt{-2}\, \theta^+\theta^-
        (\overline{\theta}^+\overline{\lambda}_{a,+}
          + \overline{\theta}^-\overline{\lambda}_{a,-})
     - \sqrt{-2}\, \overline{\theta}^+\overline{\theta}^-
        (\theta^-\lambda_{a,-} + \theta^+\lambda_{a,+})
     -2\,\theta^+\theta^-\overline{\theta}^+\overline{\theta}^-\,
         D_a\,,\\[.6ex]
   & & \hspace{-2em}a=1,\,\ldots,\, n-d\,,
  \end{eqnarray*}
 {\normalsize gauge}} 
 group $U(1)^{n-d}$ (parameterized by $(t_1,\,\ldots,\,t_{n-d})$),
 and a $U(1)^{n-d}$-action on $\Phi_i$ by
 $$
  \Phi_i\;\longrightarrow\;
          \left(\prod_{a=1}^{n-d}\,t_a^{Q_{i,a}}\right)\,\Phi_i\,.
 $$
Define ([Wi1$\,$: Eq.(2.16)]), (cf.\ curvatures)

 \vspace{-1ex}
 {\footnotesize
 \begin{eqnarray*}
  \lefteqn{
   \Sigma_a\; :=\; \frac{1}{\sqrt{2}}\overline{D}_+D_-V_a } \\[.6ex]
     & &
     =\;\sigma_a-\sqrt{-2}\theta^+\overline{\lambda}_{a,+}
      -\sqrt{-2}\,\overline{\theta}^-\lambda_{a,-}
      + \sqrt{2}\theta^+\overline{\theta}^-(D_a-\sqrt{-1}v_{a,01}) \\[.6ex]
     & & \hspace{2em}
      -\sqrt{-1}\,\overline{\theta}^-\theta^-
                  (\partial_0-\partial_1)\sigma_a
      -\sqrt{-1}\theta^+\overline{\theta}^+
                  (\partial_0+\partial_1)\sigma_a     \\[.6ex]
     & & \hspace{2em}
      +\sqrt{2}\overline{\theta}^-\theta^+\theta^-(\partial_0-\partial_1)
         \overline{\lambda}_{a,+}
      +\sqrt{2}\theta^+\overline{\theta}^-\overline{\theta}^+
         (\partial_0+\partial_1)\lambda_{a,-}
      -\theta^+\overline{\theta}^-\theta^-\overline{\theta}^+
         (\partial_0^2-\partial_1^2)\sigma_a\,,
  \end{eqnarray*}
 {\normalsize where}} 
 $v_{a,01}=\partial_0 v_{a,1}-\partial_1 v_{a,0}\,$.
The associated {\it gauged linear sigma model} is a $2$-dimensional
 supersymmetric quantum field theory (SQFT) with action
 $$
  L\; =\;    L_{\,\mbox{\scriptsize kinetic}}\,
         +\, L_{\,W}\,
         +\, L_{\,\mbox{\scriptsize gauge}}\,
         +\, L_{D,\theta}\,,
 $$
 where
 $$
   L_{\rm\, kinetic}\; =\; \int d^2y\,d^4\theta\,
   \sum_{i=1}^n
    \overline{\Phi}_i\,\exp\left[2\sum_{a=1}^{n-d}
                                Q_{i,a}V_a\right]\,\Phi_i\,,
 $$
 $$
  L_W\;
   =\; -\int d^2y\, d\theta^+\, d\theta^-\,
        W(\Phi_i)|_{\overline{\theta}^+=\overline{\theta}^{\,-}=0}
       - (\,\mbox{\rm Hermitian conjugate}\,)\,,
 $$
 $$
  L_{\,\mbox{\scriptsize gauge}}\;
   =\; -\sum_{a=1}^{n-d} \frac{1}{4\,e_a^2}\,\int d^2y\,d^4\theta\,
          \overline{\Sigma}_a\,\Sigma_a\,,
 $$
 and
 $$
  L_{D,\theta}\;
   =\;\sum_{a=1}^{n-d}\int d^2y\,
            (-r_a D_a+\frac{\theta_a}{2\pi}v_{a,01})\; \hspace{19em}
 $$
 \vspace{-1ex}
 $$
   \hspace{1.6em}=\;\frac{\sqrt{-1}\,t_a}{2\sqrt{2}}
      \int d^2y\, d\theta^+\, d\overline{\theta}^{\,-}
        \Sigma|_{\theta^-=\overline{\theta}^+=0}\,
      -\,\frac{\sqrt{-1}\,\overline{t_a}}{2\sqrt{2}}
       \int d^2y\, d\theta^-\, d\overline{\theta}^{+}
          \overline{\Sigma}|_{\theta^+=\overline{\theta}^{\,-}=0}
 $$
 with
 $$
  t_a\;=\;\sqrt{-1}\,r_a+\frac{\theta_a}{2\pi}\,.
 $$
The real-valued numerical parameters $e_a$, $r_a$, and $\theta_a$
 are called the {\it coupling constants} of the theory.

Performing the Fermionic integrations $\int d^4\theta$,
 $\int d\theta^+d\theta^-$, and
 $\int d\overline{\theta}^+d\overline{\theta}^-$ renders
 $L$ a complicated expression in terms of component fields on
 $\Phi_i$ and $V_a$ ([Wi1$\,$: Eq.(2.19), Eq.(2.21), and Eq.(2.23)]).
The $d=2, N=(2,2)$ SUSY algebra without R-symmetry generators can be
 realized as an algebra of derivations acting on $V_a$ and on $\Phi_i$
 with gauge transformations taken into account while the R-symmetry
 acts on fields via global abelian transformations on the fields.
In particular, the SUSY transformations in terms of component
 fields are given in [Wi1$\,$: Eq.(2.12) and Eq.(2.14)]
 (Cf.\ Example A.2 (1), [W-B$\,$ Chapters III-VII], and
       [We$\,$: Sec.\ 2.7.8]).
The action $L$ is invariant under these transformations
 (and hence supersymmetric).

\bigskip

\noindent $\bullet$
{\bf Wick rotation.}
Field theories on Riemannian manifolds behave better
 than those on Lorentzian manifolds.
A {\it Wick rotation} is meant to be an analytic continuation
 between theories in the two categories (e.g.\ [P-S]).
Some of its geometry is studied in [Liu].
In the current case of flat space-times, such an analytic
 continuation is realized by setting $y^0=-iy^2$ and taking
 $(y^1,y^2)$ as the coordinates of the Wick rotated $d=2$ space-time.
The latter has the Euclidean metric
 $-(dy^0)^2+(dy^1)^2=(dy^2)^2+(dy^1)^2$ and
the tangent bundle group $\SO(1,1)$ now becomes $\SO(2)$.

\bigskip

\noindent $\bullet$
{\bf The A-twist and the B-twist.}
([Wi2], [Wi1], [DEFJKMMW$\,$: vol.\ 2, Witten's lecture, Sec.\ 14.3],
  and [F-S$\,$: Chapter 7].)
 Consider the two different twisted embeddings of the $d=2$ rotation
  algebra, generated by $L$, into the $d=2,\, N=(2,2)$ SUSY algebra$\,$:
 $$
  \begin{array}{ll}
   \mbox{A-twist}\,:
     & L\;\mapsto\; L^A:= L-\frac{1}{2}J_L + \frac{1}{2}J_R \\[.6ex]
   \mbox{B-twist}\,:
     & L\;\mapsto\; L^B:= L+\frac{1}{2}J_L + \frac{1}{2}J_R\,.
  \end{array}
 $$
Since both $J_L$ and $J_R$ commute with all the SUSY algebra generators
 except $Q_{\pm}$ and $\overline{Q}_{\pm}$, the commutation relation of
 $L^A$ and $L^B$ with SUSY algebra generators are the same as those
 for $L$ except the following ones$\,$:
 $$
  \begin{array}{llll}
   [L^A,Q_+]\;=\; Q_+\,,
     & [L^A,Q_-]\;=\; 0\,,
     & [L^A,\overline{Q}_+]\;=\; 0\,,
     & [L^A,\overline{Q}_-]\;=\; -\,\overline{Q}_-\,,  \\[.6ex]
   [L^B,Q_+]\;=\; Q_+\,,
     & [L^B,Q_-]\;=\; -\,Q_-\,,
     & [L^B,\overline{Q}_+]\;=\; 0\,,
     & [L^B,\overline{Q}_-]\;=\; 0\,.
  \end{array}
 $$
This implies that all the SUSY generators are now of integral
 spin with respect to either of the twisted tangent bundle groups.

From the commutation relations of the Lorentz generator with SUSY
 generators, the supermanifolds associated to the Wick-rotated
 $d=2,\, N=(2,2)$ SUSY algebra are
 $$
   X\; =\; \left(\,C,\,\Sym^{\bullet}
         \prod \left( (K_C^{\frac{1}{2}})^{\,\oplus 2}
           \oplus (K_C^{-\frac{1}{2}})^{\,\oplus 2}\right)\right)\,,
 $$
 where $C$ is a Riemann surface and $K_C$ is the canonical line bundle
 of $C$
while the supermanifolds associated to either A-twisted or B-twisted
 SUSY algebra are
 $$
  X^{\,\mbox{\scriptsize twist}}\;
   =\; \left(C,\, \Sym^{\bullet}\prod\left(
    {\cal O}_C^{\;\oplus 2}\oplus K_C\oplus K_C^{-1}\right)\right)\,.
 $$
When $C$ is the complex plane, a cylinder, or an elliptic curve,
 $K_C^{\pm\frac{1}{2}}\sim {\cal O}_C$ admit nontrivial global sections.
Thus both SQFT and its twists can be built on such $C$.
For general $C$, $K_C^{\pm\frac{1}{2}}$ have no global sections except
 the zero-section and, hence, only twisted SQFT can be defined on $C$.

\bigskip

\noindent $\bullet$
{\bf Phase structure.}
([Wi1] and [M-P] for GLSM;
  [Al], [Arg], [Po1], [R-S-Z], and [W-K] for general field-theoretical
  aspects.)
When the coupling constants $(e_a,r_a,\theta_a)_a$ in the action $L$
 of the gauged linear model given earlier vary, the nature of the field
 theories may also vary. Thus $(e_a, r_a,\theta_a)_a$ can be thought of
 as the coordinates for a space ${\cal M}_{GLSM}$ that parameterizes
 a family of $d=2$ field theories. ${\cal M}_{GLSM}$ is called the
 ({\it Wilson's}) {\it theory space} of the model.
Quantities (e.g.\ $2$-point functions) of the field theories may be
 turned into defining geometric data (e.g.\ Zamolodchikov metric)
 on ${\cal M}_{GLSM}$.
There exists a stratification of ${\cal M}_{GLSM}$ according to the
 nature of the field theory a point in ${\cal M}_{GLSM}$ parameterizes.
Each stratum of this stratification is called a {\it phase} of the GSLM.
In general, for quantum considerations of the theory, cutoffs
(e.g.\ of energy) may have to be introduced. These cutoff parameters
 may also be added in to enlarge the theory space.
(Cf.\ The recent work of Borchard [Bo] enlarges Wilson's theory 
       space by adding also the space of 
       {\it renormalization prescriptions}.
      His work should be important to understanding the quantum phase
       structure on the theory space.)

\bigskip

\noindent $\bullet$
{\bf The moduli space of the A-twisted theory in the geometric phase.}
Either twist breaks half of the $4$ supersymmetries in general.
The resulting $d=2$ SQFT has the same expression as the action $L$
 but each of the component fields in the superfield involved lives
 in a new bundle determined by the twisted spin discussed in
 Item (The A-twist and the B-twist) above.
For the $A$-twist, the remaining supersymmetries of the twisted gauged
 linear sigma model are generated by $Q_-$ and $\overline{Q}_+\,$.
 These are the SUSY generators that are of A-twisted spin $0$.
Recall the realization of SUSY algebra with R-symmetry generators
 removed as an algebra of derivations acting on fields.
A field configuration that is annihilated by both $Q_-$ and
 $\overline{Q}_+$ is called a {\it supersymmetric field configuration}
 of the twisted gauged linear sigma model.
When the gauge coupling constants $e_a$ are all set equal to some $e$
 and the superpotential $W$ is set to zero, the bosonic part of
 SUSY configurations for the A-twisted theory are given explicitly
 by the solutions to the following system of equations
 ([Wi1$\,$: Eq.(3.33), Eq.(3.34), Eq.(3.35)] and
  [M-P$\,$: Eq.(3.54 a-d)])
 $$
  \begin{array}{rclll}
   d\sigma_a    & =  & 0\,,  & a=1\,,\,\ldots\,, n-d\,, & \\
   \sum_{\,a=1}^{n-d}\,Q^a_i\,\sigma_a\,\phi_i  & = & 0\,,
                &  i=1\,,\,\ldots\,, n\,, & \\
   D_{\overline{z}}\,\phi_i  & = & 0\,,  &  i=1\,,\,\ldots\,, n\,,
                             & \hspace{4em}(\,\ast 1\,)\\
   D_a + v_{a,12}   & = & 0\,, &  a=1\,,\,\ldots\,, n-d\,,
                             & \hspace{4em}(\,\ast 2\,)
  \end{array}
 $$
 with $D_{\overline{z}}$ a covariant derivative constructed from
  the $U(1)^{n-d}$ gauge connection $(v_{a,1},v_{a,2})_a$ and
 $$
  D_a\;=\;-e^2\,\left(\,\sum_{i=1}^n Q_i^a\,|\phi_i|^2-r_a\,\right)
 $$
 from the equation of motion for $D_a$.

When $(r_a)_a$ is in the geometric phase, for which the solution set
 to the subsystem $\{\,\mbox{Eq.($\ast 1$), Eq.($\ast 2$)}\,\}$ is
 non-empty, the only solution for $\sigma_a$ is $\sigma_a=0$ for all $a$.
Following the study in [Brad], [B-D], and [GP] on vortex-type equations,
 Witten and Morrison-Plesser thus conclude that the moduli space of
 the A-twisted theory for $(r_a)_a$ in this phase is given by

 \vspace{-2ex}
 {\small
 $$
  \begin{array}{rcl}
    \coprod_{\vec{d}}\,{\cal M}_{\vec{d}}
     & = & \left.\left\{\,\mbox{\parbox{22ex}{common solutions to
            \newline Eq.\ $(\ast 1)$ and Eq.\ $(\ast 2)$}}\,\right\}
      \right/\!\mbox{\raisebox{-1ex}{$\left\{\,\mbox{\parbox{43ex}{unitary
          abelian gauge transformations \newline
          and global complex abelian transformations}}\,\right\}$}}
                                                            \\[4ex]
     & = & \left.\left\{\,\mbox{\parbox{32ex}{solutions to
                  Eq.\ $(\ast 1)$ that satisfy \newline
                  appropriate stability condition}}\,\right\}
      \right/\!\mbox{\raisebox{-.4ex}{$\left\{\,\mbox{\parbox{23ex}{complex
               abelian gauge\newline transformations}}\,\right\}$}}\\[3ex]
     & = &  \mbox{the toric variety}\hspace{2ex}
            \coprod_{\vec{d}}\,(\,Y_{\vec{d}} - F_{\vec{d}}\,)\,/\,G\,,
  \end{array}
 $$
 {\normalsize where}} 
 $Y_{\vec{d}}$, $F_{\vec{d}}$, and $G$ are explained in
 Explanation/Fact 2.1.2  
 and Sec.\ 3 in terms of toric geometry ([M-P$\,$: Sec.\ 3.1]).
(Cf.\ See also [Fr$\,$: Lecture 4 and Lecture 5] on
      the {\it moduli space of vacua} of a SQFT.)

As already mentioned in [M-P$\,$: Sec.\ 3.7], the above construction,
 in particular the moduli space $\coprod_{\vec{d}}{\cal M}_{\vec{d}}\,$,
 has a generalization to higher genus Riemann surfaces as well,
 following [Cox2].
The main theme of this paper is the study of this generalization of
 $\coprod_{\vec{d}}{\cal M}_{\vec{d}}$ to higher genus.

\bigskip

\newpage

\end{document}